\documentclass[
	english
]{scrartcl}

\usepackage[T1]{fontenc}
\usepackage[utf8]{inputenc}

\usepackage{amsmath,amsfonts,amsthm,amssymb}
\usepackage[colorlinks,citecolor=blue]{hyperref}
\usepackage{doi,natbib}
\usepackage{natbib}
\usepackage{cancel}
\usepackage{relsize}

\usepackage{changes}

\usepackage{enumitem}
\setlist[enumerate]{label=(\alph*)}

\usepackage[figure]{hypcap}
\usepackage{graphicx}
\graphicspath{{./img/}}
\usepackage{subcaption} 

\usepackage{preamble}

\usepackage{marginnote}

\usepackage[capitalise,nameinlink]{cleveref}
\allowdisplaybreaks

\numberwithin{equation}{section}

\newcommand\norm[1]{\left\Vert#1\right\Vert}
\newcommand\nnorm[1]{\Vert#1\Vert}

\newcommand\innerprod[2]{\left\langle #1, #2\right\rangle}

\newcommand\N{\mathbb{N}}
\newcommand\R{\mathbb{R}}

\newcommand\SSS{\mathcal S}
\newcommand\PPP{\mathbf P}
\newcommand\QQQ{\mathbf Q}
\newcommand\OOO{\mathbf O}
\newcommand\LLL{\mathbf\Lambda}

\newcommand\tto{\rightrightarrows}

\newcommand{\spa}{\operatorname{span}}

\newcommand{\dist}{\operatorname{dist}}

\newcommand{\dom}{\operatorname{dom}}
\renewcommand{\Im}{\operatorname{Im}}
\newcommand{\gph}{\operatorname{gph}}
\newcommand{\epi}{\operatorname{epi}}

\newcommand{\trace}{\operatorname{trace}}

\newcommand{\xb}{\bar x}

\newcommand{\yb}{\bar y}

\DeclareMathAlphabet{\mathpzc}{OT1}{pzc}{m}{it}
\newcommand\oo{\mathpzc{o}}

\newtheorem{theorem}{Theorem}[section]
\newtheorem{lemma}[theorem]{Lemma}
\newtheorem{proposition}[theorem]{Proposition}

\newtheorem{corollary}[theorem]{Corollary}
\newtheorem{remark}[theorem]{Remark}
\newtheorem{definition}[theorem]{Definition}
\newtheorem{example}[theorem]{Example}

\definecolor{mygreen}{rgb}{0.0,0.7,0.0}
\definecolor{mybrown}{rgb}{0.5,0.5,0.0}

\begin{document}

\title{On the directional asymptotic approach in optimization theory\\
	Part A: approximate, M-, and mixed-order stationarity
	}
\author{%
	Mat\'{u}\v{s} Benko%
	\footnote{%
		University of Vienna,
		Applied Mathematics and Optimization,
		1090 Vienna,
		Austria,
		\email{matus.benko@univie.ac.at},
		\url{https://www.mat.univie.ac.at/\~rabot/group.html}
		}
	\and
	Patrick Mehlitz%
	\footnote{%
		Brandenburg University of Technology Cottbus--Senftenberg,
		Institute of Mathematics,
		03046 Cottbus,
		Germany,
		\email{mehlitz@b-tu.de},
		\url{https://www.b-tu.de/fg-optimale-steuerung/team/dr-patrick-mehlitz},
		\orcid{0000-0002-9355-850X}%
		}
	}

\publishers{}
\maketitle

\begin{abstract}
	 We show that, for a fixed order $\gamma\geq 1$, each local minimizer of
	 a rather general nonsmooth optimization problem in Euclidean spaces
	 is either M-stationary in the classical sense (corresponding
	 to stationarity of order $1$), satisfies
	 stationarity conditions in terms of a coderivative construction of
	 order $\gamma$, or is approximately stationary with respect
	 to a critical direction as well as $\gamma$ in a certain sense. 
	 By ruling out the latter case 
	 with a constraint qualification not stronger than directional metric
	 subregularity, we end up with new necessary optimality conditions comprising
	 a mixture of limiting variational tools of order $1$ and $\gamma$.
	 These abstract findings are carved out for the broad class of geometric
	 constraints. As a byproduct, we obtain new constraint qualifications
	 ensuring M-stationarity of local minimizers.
	 The paper closes by illustrating these results in the context
	 of standard nonlinear, complementarity-constrained, and nonlinear
	 semidefinite programming.
\end{abstract}

\begin{keywords}	
	Approximate stationarity, Complementarity-constrained programming, 
	Directional limiting variational calculus, M-stationarity, 
	Pseudo-coderivatives, Semidefinite programming
\end{keywords}

\begin{msc}	
	\mscLink{49J52}, \mscLink{49J53}, \mscLink{49K27}, \mscLink{90C22},
	\mscLink{90C30}, \mscLink{90C33}
\end{msc}

\section{Introduction}\label{sec:introduction}

In order to identify local minimizers of optimization problems analytically
or numerically, it is
desirable that such points satisfy applicable necessary optimality conditions.
Typically, under validity of a constraint qualification, first-order
necessary optimality conditions of abstract Karush--Kuhn--Tucker (KKT)-type
hold at local minimizers. Here,
\emph{first-order} refers to the fact that first-order tools of (generalized)
differentiation are used to describe the variation of all involved data functions.
In case where the celebrated tools of limiting variational analysis are exploited,
one speaks of so-called Mordukhovich (or, briefly, M-) stationarity,
see \cite{Mordukhovich2018}.
In the absence of constraint qualifications, i.e., in a \emph{degenerate} situation,
local minimizers still satisfy a Fritz--John (FJ)-type first-order necessary
optimality condition which allows for a potentially vanishing multiplier
associated with the generalized derivative of the objective function.
Since such a condition allows to discard the objective function, it might be
too weak in practically relevant scenarios. 
In \cite{KrugerMehlitz2021,Mehlitz2020a}, the authors
introduce an \emph{approximate} (sometimes referred to as asymptotic) 
necessary optimality condition of first order which
is an intermediate concept between the aforementioned KKT- and FJ-type 
M-stationarity conditions
for rather general optimization problems. On the one
hand, approximate stationarity holds at local minimizers without a constraint
qualification and it possesses notable relevance in numerical optimization when
convergence properties of solution algorithms are investigated. 
On the other hand, it is not easy to check validity of this condition in practice
since it is stated in terms of sequences whose iterates satisfy certain stationarity
conditions up to a residual which tends to zero.
We would like to note that approximate stationarity corresponds to the so-called
AKKT conditions in standard nonlinear programming, see e.g.\
\cite{AndreaniMartinezSvaiter2010,AndreaniHaeserMartinez2011}, 
which have been extended to bilevel, disjunctive, conic, 
and even infinite-dimensional optimization,
recently, see e.g.\ \cite{AndreaniHaeserSecchinSilva2019,AndreaniGomezHaeserMitoRamos2021,BoergensKanzowMehlitzWachsmuth2019,KanzowSteckWachsmuth2018,KanzowRaharjaSchwartz2021,Mehlitz2022,Ramos2019}
and the references therein.

The aim of this paper is to apply the \emph{directional} approach
to limiting variational analysis, see e.g.\ \cite{BenkoGfrererOutrata2019}, 
in order to enrich the approximate
stationarity conditions from \cite{KrugerMehlitz2021,Mehlitz2020a}
with the aid of directional information.
Noting that the directional tools of variational analysis were successfully applied 
to find refined M-stationarity-type optimality conditions and mild constraint qualifications 
for diverse problems in optimization theory, see e.g.\
\cite{BaiYeZhang2019,BaiYe2021,BenkoCervinkaHoheisel2019,Gfrerer2013,Gfrerer2014,GfrererKlatte2016,GfrererYeZhou2022}
and the references therein, this seems to be a desirable goal.

Based on the tools of limiting variational analysis, we show 
in \cref{thm:higher_order_directional_asymptotic_stationarity} that local minimizers
of rather general optimization problems in Euclidean spaces are either
M-stationary, satisfy a stationarity condition combining the limiting subdifferential
of the objective function and a coderivative-like tool associated with the
constraints of some arbitrary order $\gamma\geq 1$, 
a so-called \emph{pseudo-coderivative}, see \cite{Gfrerer2014a}, or come along with an
approximate stationarity condition depending on a critical direction as well as
the order $\gamma$ where the involved sequence of multipliers is diverging. 
Even for $\gamma:=1$, this enhances the findings from
\cite{KrugerMehlitz2021,Mehlitz2020a}. Furthermore, this result opens a new way on
how to come up with applicable necessary optimality conditions for the original problem,
namely, by ruling out the irregular situation of approximate stationarity which can
be done in the presence of so-called \emph{metric pseudo-subregularity} of order $\gamma$,
see \cite{Gfrerer2014a} again.
In case $\gamma:=1$, we end up with M-stationarity, and metric pseudo-subregularity reduces
to metric subregularity, i.e., we obtain results related to \cite{Gfrerer2013}.
For $\gamma>1$, this procedure leads to a mixed-order stationarity condition 
involving the pseudo-coderivative of order $\gamma$, and metric pseudo-subregularity
is weaker than metric subregularity.
If $\gamma:=2$ and so-called geometric constraints, induced by a twice continuously
differentiable mapping $g$ as well as a closed set $D$, are investigated, 
this pseudo-coderivative
can be estimated from above in terms of initial problem data, i.e., 
in terms of (first- and second-order) derivatives associated with $g$ 
as well as tangent and normal cones to $D$, under mild conditions, 
see \cref{sec:variational_analysis_constraint_mapping}. The associated
mixed-order necessary optimality conditions and qualification conditions are worked
out in \cref{sec:constraint_mappings}. 
Let us note that necessary optimality conditions for 
degenerate optimization problems which comprise 
a mixture of first- and second-order information at the same time can be found in
\cite{Avakov1989,AvakovArutunovIzmailov2007,Gfrerer2007,Gfrerer2014a}.
All these results comprise a multiplier associated with the objective function
which is allowed to be zero. Our approach yields validity of such necessary
optimality conditions with leading multiplier $1$ under mild qualification conditions
which can be checked in terms of initial problem data.
Furthermore, we show how our results can
be used to find new constraint qualifications guaranteeing M-stationarity of
local minimizers.
In Part B of this paper, the same conditions are derived from a different angle
and are shown to be not stronger than the prominent 
\emph{First-} and \emph{Second-Order Sufficient Condition for Metric Subregularity},
see \cite{GfrererKlatte2016}.
Exemplary, we illustrate some of these findings in \cref{sec:applications} 
by means of standard nonlinear, 
complementarity-constrained, and nonlinear semidefinite programs. 
Notably, the latter setting is non-polyhedral due to the appearance of the cone
of all positive semidefinite matrices which possesses curvature.

The remaining parts of the paper are organized as follows:
In \cref{sec:notation}, we document the notation used in this paper and provide
some preliminary results from variational analysis and generalized differentiation.
Moreover, in \cref{sec:generalized_differentiation}, we introduce and investigate two
new generalized derivatives for set-valued mappings, a primal construction
called \emph{graphical subderivative} and a dual construction called
\emph{directional limiting pseudo-coderivative} which slightly differs from the tool
used in \cite{Gfrerer2014a}.
\Cref{sec:variational_analysis_constraint_mapping} is dedicated to the study of the generalized 
derivatives of so-called constraint mappings which encode most of the constraints
which appear in optimization theory. Particularly, we present upper estimates
for the directional limiting pseudo-coderivative of such mappings.
In \cref{sec:directional_asymptotic_tools}, we first derive the aforementioned 
stationarity condition for abstract nonsmooth optimization problems which holds 
in the absence of constraint qualifications and shows that local minimizers are
either M-stationary, stationary with respect to (w.r.t.) 
the pseudo-coderivative, or approximately
stationary. Furthermore, we demonstrate that the latter case can be ruled out in
the presence of directional metric pseudo-subregularity which, thus, leads to
applicable optimality conditions provided the pseudo-coderivative can be computed
or estimated. Since this is possible for constraint maps, we specify our findings
for these special mappings in \cref{sec:constraint_mappings}. Additionally, this
approach yields new constraint qualifications guaranteeing M-stationarity
of local minimizers. 
In order to illustrate our findings for particularly interesting problem classes,
some consequences of the derived theory in the setting of standard nonlinear,
complementarity-constrained, and nonlinear semidefinite programming are
presented in \cref{sec:applications}.
The paper closes with some concluding remarks in \cref{sec:conclusions}.

\section{Notation and preliminaries}\label{sec:notation}

We rely on standard notation taken from
\cite{AubinFrankowska2009,BonnansShapiro2000,RockafellarWets1998,Mordukhovich2018}.

\subsection{Basic notation}

Throughout the paper, $\mathbb X$ and $\mathbb Y$ denote Euclidean spaces, 
i.e., finite-dimensional
Hilbert spaces. For simplicity, the associated inner product will be represented by
$\innerprod{\cdot}{\cdot}$ since the underlying space will be clear from the context.
The norm induced by the inner product is denoted by $\norm{\cdot}$.
The unit sphere in $\mathbb X$ will be represented by $\mathbb S_{\mathbb X}$ throughout the paper.
Furthermore, for $\varepsilon>0$ and $\bar x\in\mathbb X$,
$\mathbb B_\varepsilon(\bar x):=\{x\in\mathbb X\,|\,\norm{x-\bar x}\leq\varepsilon\}$
is the closed $\varepsilon$-ball around $\bar x$.
For a nonempty set $Q\subset\mathbb X$, 
$Q^\circ:=\{\eta\in\mathbb X\,|\,\forall x\in Q\colon\,\innerprod{\eta}{x}\leq 0\}$
is the so-called polar cone of $Q$ which is always a closed, convex cone.
Furthermore, for some $\bar x\in\mathbb X$, 
$[\bar x]^\perp:=\{\eta\in\mathbb X\,|\,\innerprod{\eta}{\bar x}=0\}$ is the annihilator
of $\bar x$. 
For simplicity of notation, we use $\bar x+Q:=Q+\bar x:=\{x+\bar x\in\mathbb X\,|\,x\in Q\}$.
For a given linear operator $A\colon\mathbb X\to\mathbb Y$, 
$A^*\colon\mathbb Y\to\mathbb X$ is used to denote its adjoint.

Let $g\colon\mathbb X\to\mathbb Y$ be a continuously differentiable mapping.
We use $\nabla g(\bar x)\colon\mathbb X\to\mathbb Y$ in order to denote the derivative
of $g$ at $\bar x\in\mathbb X$. Note that $\nabla g(\bar x)$ is a linear operator.
Let us emphasize that, in the special case $\mathbb Y:=\R$, $\nabla g(\bar x)$
does not coincide with the standard gradient which would correspond to
$\nabla g(\bar x)^*1$.
For twice continuously differentiable $g$ and a scalar $\lambda\in\mathbb Y$,
we set $\langle \lambda,g\rangle(x):=\innerprod{\lambda}{ g(x) }$ 
for each $x\in\mathbb X$ in order
to denote the associated scalarization mapping $\innerprod{\lambda}{g}\colon\mathbb X\to\R$.
By $\nabla\innerprod{\lambda}{g}(\bar x)$ and $\nabla^2\innerprod{\lambda}{g}(\bar x)$
we represent the first- and second-order derivatives of this map at $\bar x\in\mathbb X$ (w.r.t.\ the variable which enters $g$).
Furthermore, for $u\in\mathbb X$, we make use of
\[
	\nabla^2g(\bar x)[u,u]
	:=
	\sum_{i=1}^m \innerprod{u}{\nabla^2\innerprod{e_i}{g}(\bar x)(u)}\,e_i
\]
for brevity where 
$m\in\N$ is the dimension of $\mathbb Y$ and
$e_1,\ldots,e_m\in\mathbb Y$ denote the $m$ canonical unit vectors of $\mathbb Y$.
In case $\mathbb Y:=\R$, the second-order derivative 
$\nabla ^2g(\bar x)\colon\mathbb X\times\mathbb X\to\R$ is a bilinear mapping, and
for each $u\in\mathbb X$, we identify $\nabla^2g(\bar x)u$ with an element of $\mathbb X$.

\subsection{Fundamentals of variational analysis}

Let us fix a closed set $Q\subset\mathbb X$ and some point $x\in Q$.
We use
\begin{align*}
	\mathcal T_Q(x)
	:=
	\left\{d\in\mathbb X\,\middle|\,
		\begin{aligned}
			&\exists\{d_k\}_{k\in\N}\subset\mathbb X,\,\exists\{t_k\}_{k\in\N}\subset\R_+\colon\\
			&\qquad d_k\to d,\,t_k\searrow 0,\,x+t_kd_k\in Q\,\forall k\in\N
		\end{aligned}
	\right\}
\end{align*}
to denote the (Bouligand) tangent cone to $Q$ at $x$.
Furthermore, we make use of
\begin{align*}
	\widehat{\mathcal N}_Q(x)
	&:=
	\left\{\eta\in\mathbb X\,|\,\forall x'\in Q\colon\,\innerprod{\eta}{x'-x}\leq\oo(\nnorm{x'-x})\right\},\\
	\mathcal N_Q(x)
	&:=
	\left\{\eta\in\mathbb X\,\middle|\,
		\begin{aligned}
			&\exists\{x_k\}_{k\in\N}\subset Q,\,\exists\{\eta_k\}_{k\in\N}\subset\mathbb X\colon\\
			&\qquad x_k\to x,\,\eta_k\to\eta,\,\eta_k\in\widehat{\mathcal N}_Q(x_k)\,\forall k\in\N
		\end{aligned}
	\right\},
\end{align*}
the regular (or Fr\'{e}chet) and limiting (or Mordukhovich) normal cone to $Q$ at $x$.
Observe that both of these normal cones coincide with the standard normal cone of convex
analysis as soon as $Q$ is convex.
For $\tilde x\notin Q$, we set $\mathcal T_Q(\tilde x):=\varnothing$ and
$\widehat{\mathcal N}_Q(\tilde x)=\mathcal N_Q(\tilde x):=\varnothing$.
Finally, for some $d\in\mathbb X$, we use
\[
	\mathcal N_Q(x;d)
	:=
	\left\{
		\eta\in\mathbb X\,\middle|\,
			\begin{aligned}
				&\exists\{d_k\}_{k\in\N}\subset\mathbb X,\,\exists\{t_k\}_{k\in\N}\subset\R_+,\,
				\exists\{\eta_k\}_{k\in\N}\subset\mathbb X\colon\\
				&\qquad d_k\to d,\,t_k\searrow 0,\,\eta_k\to\eta,\,\eta_k\in\widehat{\mathcal N}_Q(x+t_kd_k)\,\forall k\in\N
			\end{aligned}
	\right\}
\]
in order to represent the directional limiting normal cone to $Q$ at $x$ in direction $d$.
Note that this set is empty if $d$ does not belong to $\mathcal T_Q(x)$. 
If $Q$ is convex, we have $\mathcal N_Q(x;d)=\mathcal N_Q(x)\cap [d]^\perp$.

In this paper, the concept of polyhedrality will be of essential importance.
Let us recall that a set $Q\subset\R^m$ will be called polyhedral if it is
the union of finitely many convex polyhedral sets. Similarly, it is referred to
as locally polyhedral around $x\in Q$ whenever 
$Q\cap\{z\in\R^m\,|\,\forall i\in\{1,\ldots,m\}\colon\,|z_i-x_i|\leq\varepsilon\}$
is polyhedral for some $\varepsilon>0$.
The following lemma provides some basic properties of polyhedral sets.
Statement~\ref{item:exactness_tangential_approximation}, 
proven in \cite[Proposition 8.24]{Ioffe2017},
is an extension of the exactness of tangential approximations 
of convex polyhedral sets from \cite[Exercise~6.47]{RockafellarWets1998}.
Thus, we refer to this property of (locally) polyhedral sets as 
\emph{exactness of tangential approximations} again.
The equality in statement~\ref{item:normal_cones_to_polyhedral_sets} 
follows from \cite[Proposition 2.11]{BenkoGfrererYeZhangZhou2022}
and the rest is straightforward,
see \cite[Lemma~2.1]{Gfrerer2014} as well.
\begin{lemma}\label{lem:some_properties_of_polyhedral_sets}
	Let $Q\subset\mathbb \R^n$ be a closed set which is locally polyhedral
	around some fixed $x\in Q$.
	Then the following statements hold.
	\begin{enumerate}
	\item\label{item:exactness_tangential_approximation}
		There exists a neighborhood $U\subset\R^n$ of $x$ 
		such that $(x+\mathcal T_Q(x))\cap U=Q\cap U$.
	\item\label{item:normal_cones_to_polyhedral_sets} 
		For arbitrary $w\in\R^n$, we have
		\begin{equation}\label{eq:normals_to_polyhedral_sets}
			\mathcal N_Q(x;w)
			=
			\mathcal N_{\mathcal T_Q(x)}(w)
			\subset
			\mathcal N_{Q}(x)\cap [w]^\perp.
		\end{equation}
		If $Q$ is, additionally, convex, and $w\in\mathcal T_Q(x)$, 
		then the final inclusion holds as an equality.
	\end{enumerate}
\end{lemma}

It is well known that the regular and limiting normal cone enjoy an exact product rule
which is not true for the tangent cone in general.
However, the following lemma shows that such a product rule also holds for tangents
as soon as polyhedral sets are under consideration.
Its proof is straightforward and, hence, omitted.
\begin{lemma}\label{lem:product_rule_tangents_polyhedral_sets}
\leavevmode
	\begin{enumerate}
		\item For closed sets $P\subset\mathbb X$ and $Q\subset\mathbb Y$ 
			as well as $x\in P$ and $y\in Q$,
			we have $\mathcal T_{P\times Q}(x,y)\subset\mathcal T_P(x)\times\mathcal T_Q(y)$.
		\item For closed sets $P\subset\R^n$ and $Q\subset\R^m$ as well as $x\in P$ and $y\in Q$,
			such that $P$ and $Q$ are locally polyhedral around $x$ and $y$, respectively, 
			we have $\mathcal T_{P\times Q}(x,y)=\mathcal T_P(x)\times\mathcal T_Q(y)$.
	\end{enumerate}
\end{lemma}

Let us mention that a slightly more general version of the above lemma 
can be found in \cite[Proposition~1]{GfrererYe2017}.

For a set-valued mapping $\Phi\colon\mathbb X\tto\mathbb Y$, we use
$\dom \Phi:=\{x\in\mathbb X\,|\,\Phi(x)\neq\varnothing\}$,
$\gph \Phi:=\{(x,y)\in\mathbb X\times\mathbb Y\,|\,y\in\Phi(x)\}$,
$\ker\Phi:=\{x\in\mathbb X\,|\,0\in\Phi(x)\}$, and
$\Im\Phi:=\bigcup_{x\in\mathbb X}\Phi(x)$ in order to
represent the domain, graph, kernel, and image of $\Phi$, respectively.
Furthermore, the so-called inverse mapping $\Phi^{-1}\colon\mathbb Y\tto\mathbb X$
is defined via $\gph\Phi^{-1}:=\{(y,x)\in\mathbb Y\times\mathbb X\,|\,(x,y)\in\gph\Phi\}$.

There exist numerous concepts of local \emph{regularity} or
\emph{Lipschitzness} associated with set-valued mappings.
In this paper, we are concerned with so-called
directional metric pseudo-subregularity which originates from 
\cite{Gfrerer2014a}.
\begin{definition}\label{def:metric_pseudo_subregularity}
	Fix a set-valued mapping $\Phi\colon\mathbb X\tto\mathbb Y$ which has closed graph locally
	around $(\bar x,\bar y)\in\gph\Phi$, and fix a direction $u\in\mathbb X$ and 
	a constant $\gamma\geq 1$.
	We say that $\Phi$ is \emph{metrically pseudo-subregular of order $\gamma$ in direction $u$}
	at $(\bar x,\bar y)$ if there are constants $\varepsilon>0$, $\delta>0$, and $\kappa>0$
	such that
	\[
		\forall x\in\bar x+\mathbb B_{\varepsilon,\delta}(u)\colon\quad
		\norm{x-\bar x}^{\gamma-1}\dist(x,\Phi^{-1}(\bar y))
		\leq
		\kappa\,\dist(\bar y,\Phi(x))
	\]
	where $\mathbb B_{\varepsilon,\delta}(u)
	:=\{v\in\mathbb X\,|\,\norm{\norm{v}u-\norm{u}v}\leq\delta\norm{u}\norm{v},\,
	\norm{v}\leq\varepsilon\}$
	is a so-called \emph{directional neighborhood of $u$}.
	In case where this is fulfilled for $u:=0$, we say that $\Phi$ is
	\emph{metrically pseudo-subregular of order $\gamma$} at $(\bar x,\bar y)$.
\end{definition}

For $\gamma:=1$, the above definition recovers the one of \emph{directional metric
subregularity}, see \cite{Gfrerer2013}, and if $u:=0$ is chosen additionally,
we end up with the classical notion of metric subregularity.
We also note that metric subregularity in a specified direction implies
metric pseudo-subregularity of arbitrary order $\gamma>1$ in the same direction.

\subsection{Generalized differentiation}\label{sec:generalized_differentiation}

In this section, we recall some notions from generalized differentiation and
introduce two novel derivatives for set-valued mappings.

Let us start with a lower semicontinuous function $\varphi\colon\mathbb X\to\R\cup\{\infty\}$
and some point $\bar x\in\mathbb X$ where $\varphi(\bar x)<\infty$ is valid.
Its regular (or Fr\'{e}chet) and limiting (or Mordukhovich) subdifferential at
$\bar x$ are given by
\begin{align*}
	\widehat\partial \varphi(\bar x)
	&:=
	\left\{
		\eta\in\mathbb X\,\middle|\,
		(\eta,-1)\in\widehat{\mathcal N}_{\epi\varphi}(\bar x,\varphi(\bar x))
	\right\},\\
	\partial \varphi(\bar x)
	&:=
	\left\{
		\eta\in\mathbb X\,\middle|\,
		(\eta,-1)\in\mathcal N_{\epi\varphi}(\bar x,\varphi(\bar x))
	\right\},
\end{align*}
respectively, where $\epi\varphi:=\{(x,\alpha)\in\mathbb X\times\R\,|\,\varphi(x)\leq\alpha\}$
is the epigraph of $\varphi$.
In case where $\varphi$ is continuously differentiable at $\bar x$, both sets
reduce to the singleton containing only the gradient $\nabla\varphi(\bar x)^*1$.

Below, we introduce two different graphical derivatives of a set-valued mapping. 
While the standard graphical derivative is well known from the literature, the
concept of a graphical subderivative is, to the best of our knowledge, new.
\begin{definition}\label{def:graphical_derivative}
	Let $\Phi\colon\mathbb X\tto\mathbb Y$ be a set-valued mapping possessing 
	a closed graph locally around $(\xb,\yb)\in\gph \Phi$.
	\begin{enumerate}
	\item 
		The \emph{graphical derivative} of $\Phi$ at $(\bar x,\bar y)$ is the
		mapping $D\Phi(\bar x,\bar y)\colon\mathbb X\tto\mathbb Y$ given by
		\[
			\gph D\Phi(\bar x,\bar y)=\mathcal T_{\gph\Phi}(\bar x,\bar y).
		\]
		In case where $\Phi$ is single-valued at $\bar x$, we use 
		$D\Phi(\bar x)\colon\mathbb X\tto\mathbb Y$ for brevity.
	\item 
		The \emph{graphical subderivative} of $\Phi$ at $(\bar x,\bar y)$
		is the mapping 
		$D_\textup{sub}\Phi(\bar x,\bar y)\colon \mathbb S_{\mathbb X} \tto \mathbb S_{\mathbb Y}$ 
		which assigns to every $u\in \mathbb S_{\mathbb X}$ 
		the set of all $v\in \mathbb S_{\mathbb Y}$
		such that there are sequences $\{u_k\}_{k\in\N}\subset\mathbb X$, 
		$\{v_k\}_{k\in\N}\subset\mathbb Y$, 
		and $\{t_k\}_{k\in\N},\{\tau_k\}_{k\in\N}\subset\R_+$
		which satisfy $u_k\to u$, $v_k\to v$, $t_k\searrow 0$, $\tau_k\searrow 0$, 
		$\tau_k/t_k\to\infty$, and 
		$(\bar x+t_ku_k,\bar y+\tau_kv_k)\in\gph\Phi$ 
		for all $k\in\N$.
	\end{enumerate}
\end{definition}

Let us note that for every set-valued mapping $\Phi\colon\mathbb X\tto\mathbb Y$ whose graph is
closed locally around $(\bar x,\bar y)\in\gph\Phi$, we obtain the trivial estimate
\begin{equation}\label{eq:trivial_upper_estimate_graphical_subderivative}
	\forall u\in\mathbb S_{\mathbb X}\colon\quad
	D_\textup{sub}\Phi(\bar x,\bar y)(u)\subset D\Phi(\bar x,\bar y)(0)
\end{equation}
right from the definition of these objects. 

In the course of the paper, we are mainly interested in the graphical (sub)derivative 
associated with so-called normal cone mappings. 
In the next lemma, we present some corresponding upper estimates.

\begin{lemma}\label{lem:graphical_derivatives_of_normal_cone_map}
	Let $D\subset\mathbb Y$ be a nonempty, closed, convex set such that
	the (single-valued) projection operator onto $D$, denoted by
	$\Pi_D\colon\mathbb Y\to\mathbb Y$, is directionally differentiable.
	Fix $\bar y\in D$ and $\bar y^*\in\mathcal N_D(\bar y)$.
	Then, for arbitrary $u\in\mathbb Y$, we find
	\begin{align*}
		D\mathcal N_D(\bar y,\bar y^*)(u)
		&\subset
		\{v\in\mathbb Y\,|\,\Pi_D'(\bar y+\bar y^*;u+v)=u\},
	\end{align*}
	and for $u\in\mathbb S_{\mathbb Y}$, we find
	\begin{align*}
		D_\textup{sub}\mathcal N_D(\bar y,\bar y^*)(u)
		&\subset
		\{v\in\mathbb S_{\mathbb Y}\,|\,
			\Pi_D'(\bar y+\bar y^*;v)=0, \innerprod{u}{v}\geq 0\}.
	\end{align*}
\end{lemma}
\begin{proof}
	By convexity of $D$, we have
	\[
		\forall y,y^*\in\mathbb Y\colon\quad
		y^*\in\mathcal N_D(y)
		\quad\Longleftrightarrow\quad
		\Pi_D(y+y^*)=y.
	\]
	In the remainder of the proof, we set $\tilde y:=\bar y+\bar y^*$ for brevity.
	Next, let us fix $u,v\in\mathbb Y$ as well as 
	sequences $\{u_k\}_{k\in\N},\{v_k\}_{k\in\N}\subset\mathbb Y$
	and $\{\tau_k\}_{k\in\N},\{\varepsilon_k\}_{k\in\N}\subset\R_+$
	such that $u_k\to u$, $v_k\to v$, $\tau_k\searrow 0$, 
	and $\bar y^*+\tau_kv_k\in\mathcal N_D(\bar y+\tau_k\varepsilon_k u_k)$, i.e.,
	$\Pi_D(\tilde y+\tau_k\varepsilon_ku_k+\tau_kv_k)=\bar y+\tau_k\varepsilon_ku_k$,
	for each $k\in\N$. Using $\Pi_D(\tilde y)=\bar y$, we find
	\begin{equation}\label{eq:graphical_derivative_normal_cone_map}
		\forall k\in\N\colon\quad
		\varepsilon_k u_k
		=
		\frac{\Pi_D(\tilde y+\tau_k\varepsilon_ku_k+\tau_kv_k)-\Pi_D(\tilde y)}
		{\tau_k}.
	\end{equation}
	
	In case where $v\in D\mathcal N_D(\bar y,\bar y^*)(u)$ holds, we can choose $\varepsilon_k=1$
	for each $k\in\N$, and taking the limit $k\to\infty$ in
	\eqref{eq:graphical_derivative_normal_cone_map} 
	while exploiting directional differentiability
	and Lipschitzness of $\Pi_D$ yields $\Pi'_D(\tilde y;u+v)=u$.
	This shows the first estimate.
	
	Now, assume that $v\in D_\textup{sub}\mathcal N_D(\bar y,\bar y^*)(u)$ is valid.
	Then $\varepsilon_k\searrow 0$ and $u, v \in\mathbb S_{\mathbb Y}$ can be postulated,
	and taking the limit $k\to\infty$ in \eqref{eq:graphical_derivative_normal_cone_map} 
	shows $\Pi'_D(\tilde y;v)=0$.
	By nature of the projection, we have
	\begin{align*}
		\innerprod{
			\tilde y+\tau_k\varepsilon_ku_k+\tau_kv_k
			-
			\Pi_D(\tilde y+\tau_k\varepsilon_ku_k+\tau_kv_k)
		}{
			\Pi_D(\tilde y)
			-
			\Pi_D(\tilde y+\tau_k\varepsilon_ku_k+\tau_kv_k)
		}
		\leq 
		0
	\end{align*}
	for each $k\in\N$.
	Exploiting \eqref{eq:graphical_derivative_normal_cone_map}, this is equivalent to
	\begin{align*}
		\innerprod{
			\tilde y+\tau_kv_k-\Pi_D(\tilde y)
		}{
			\Pi_D(\tilde y)
			-
			\Pi_D(\tilde y+\tau_k\varepsilon_ku_k+\tau_kv_k)
		}
		\leq 
		0
	\end{align*}
	for each $k\in\N$.
	Some rearrangements and the characterization of the projection lead to
	\begin{align*}
		&\tau_k\,
		\innerprod{v_k}{
			\Pi_D(\tilde y)
			-
			\Pi_D(\tilde y+\tau_k\varepsilon_ku_k+\tau_kv_k)
		}
		\\
		&\quad
		\leq
		\innerprod{
			\tilde y-\Pi_D(\tilde y)
		}{
			\Pi_D(\tilde y+\tau_k\varepsilon_ku_k+\tau_kv_k)
			-
			\Pi_D(\tilde y)
		}
		\leq 
		0.
	\end{align*}
	Division by $\tau_k^2\varepsilon_k$ and \eqref{eq:graphical_derivative_normal_cone_map},
	thus, give us $\innerprod{v_k}{ u_k }\geq 0$
	for each $k\in\N$, and taking the limit, we obtain $\innerprod{u}{v}\geq 0$
	which shows the second estimate.
\end{proof}

Let us note that it has been shown in \cite[Theorem~3.1, Corollary~3.1]{WuZhangZhang2014} 
that the estimate on the
graphical derivative of the normal cone mapping $\mathcal N_D$ holds as an equality 
in the situation where
$D$ is the convex cone of positive semidefinite symmetric matrices, 
and that the presented proof extends to arbitrary convex
cones as long as the associated projection operator is directionally differentiable.
This result can also be found in slightly more general form in 
\cite[Theorem~3.3]{MordukhovichOutrataRamirez2015}.
In order to make the estimates from \cref{lem:graphical_derivatives_of_normal_cone_map} 
explicit, one needs to be in
position to characterize the directional derivative of the projection onto the convex set $D$.
This is easily possible if $D$ is polyhedral, see \cite{Haraux1977} 
for this classical result, and 
\cref{rem:upper_estimate_pseudo_coderivative_constraint_maps}\,\ref{item:comparison_estimate_pseudocoderivative_polyhedral}, 
but even in
non-polyhedral situations, e.g., where $D$ is the second-order cone 
or the cone of positive semidefinite symmetric
matrices, closed formulas for this directional derivative are available in the literature, see
\cite[Lemma~2]{OutrataSun2008} and \cite[Theorem~4.7]{SunSun2002}, respectively.

In the next two results, we investigate the special situation $\mathbb Y:=\R^m$ in detail.
In case where we consider the normal cone mapping associated with polyhedral sets, 
there is no difference
between graphical derivative and graphical subderivative as the subsequent lemma shows.
\begin{lemma}\label{lem:normal_cone_map_of_polyhedral_set}
	Let $D\subset\R^m$ be a polyhedral set.
	Then $\gph\mathcal N_D$ is polyhedral as well, and 
	for arbitrary $(\bar y,\bar y^*)\in\gph\mathcal N_D$ and 
	$u, v\in\R^m \setminus \{0\}$, we have
	\[
		v \in D\mathcal N_D(\bar y,\bar y^*)(u)
		\ \iff \
		v/\norm{v} \in D_\textup{sub}\mathcal N_D(\bar y,\bar y^*)(u/\norm{u}).
	\]
\end{lemma}
\begin{proof}
	It follows from \cite[Theorem~2]{AdamCervinkaPistek2016} that there exist finitely many
	convex polyhedral sets $D_1,\ldots,D_\ell\subset\R^m$ and closed, convex, polyhedral cones
	$K_1,\ldots,K_\ell\subset\R^m$ such that $\gph\mathcal N_D=\bigcup_{i=1}^\ell D_i\times K_i$.
	Particularly, $\gph\mathcal N_D$ is polyhedral.
	
	Next, consider some nonzero $u, v\in\R^m$ with $v/\norm{v} \in D_\textup{sub}\mathcal N_D(\bar y,\bar y^*)(u/\norm{u})$. 
	Then we find $\{\tilde u_k\}_{k\in\N},\{\tilde v_k\}_{k\in\N}\subset\R^m$ and $\{\tilde t_k\}_{k\in\N},\{\tau_k\}_{k\in\N}\subset\R_+$
	such that $u_k:=\tilde u_k \norm{u}\to u$, $v_k:=\tilde v_k \norm{v} \to v$, $t_k:= \tilde t_k/\norm{u} \searrow 0$, $\tau_k\searrow 0$, $\tau_k/t_k\to\infty$, and
	$(\bar y+t_k u_k,\bar y^*+ (\tau_k/\norm{v}) v_k)\in\gph\mathcal N_D$ for all $k\in\N$.
	Thus, we can pick $j\in\{1,\ldots,\ell\}$ and a subsequence (without relabeling) such that
	$(\bar y+t_k u_k,\bar y^*+ (\tau_k/\norm{v}) v_k) \in D_j\times K_j$
	and $\tau_k/\norm{v}>t_k$ for all $k\in\N$.
	By convexity of $K_j$, we also have
	$(\bar y+t_k u_k,\bar y^*+ t_k v_k) \in D_j\times K_j$
	which shows $v\in D\mathcal N_D(\bar y,\bar y^*)(u)$.
	The converse implication can be proven in analogous fashion by multiplying the null sequence in the domain space
	with another null sequence.
\end{proof}

The next lemma shows how the graphical derivative of normal cone mappings associated with Cartesian products
of polyhedral sets can be computed.
\begin{lemma}\label{lem:product_rule_graphical_derivative}
	Fix some $\ell\in\N$. For each $i\in\{1,\ldots,\ell\}$, 
	let $D_i\subset\R^{m_i}$ for some $m_i\in\N$ be polyhedral.
	Set $D:=\prod_{i=1}^\ell D_i$, $m:=\sum_{i=1}^\ell m_i$, and $L:=\{1,\ldots,\ell\}$. 
	Then we have 
	\[
		\gph\mathcal N_D
		=
		\{((y_1,\ldots,y_\ell),(y_1^*,\ldots,y_\ell^*))\in\R^m\times\R^m\,|\,
		\forall i\in L\colon\,(y_i,y_i^*)\in\gph\mathcal N_{D_i}\},
	\]
	and for arbitrary 
	$\bar y:=(\bar y_1,\ldots,\bar y_\ell),\bar y^*:=(\bar y_1^*,\ldots,\bar y_\ell^*)\in\R^m$
	satisfying
	$(\bar y,\bar y^*)\in\gph\mathcal N_D$ and for each $u:=(u_1,\ldots,u_\ell)\in\R^m$,
	we find
	\[
		D\mathcal N_D(\bar y,\bar y^*)(u)
		=
		\{v=(v_1,\ldots,v_\ell)\in\R^m\,|\,
		\forall i\in L\colon\,v_i\in D\mathcal N_{D_i}(\bar y_i,\bar y_i^*)(u_i)\}.
	\]
\end{lemma}
\begin{proof}
	The representation of $\gph\mathcal N_D$ is a simple consequence of the product rule 
	for the computation of limiting
	normals, see e.g.\ \cite[Proposition~1.4]{Mordukhovich2018}, 
	and does not rely on the polyhedrality of the underlying
	sets. Note that this implies $\gph \mathcal N_D=A\prod_{i=1}^\ell\gph\mathcal N_{D_i}$ 
	for a suitably chosen permutation
	matrix $A\in\R^{2m\times 2m}$. 
	Particularly, $\gph \mathcal N_D$ is  polyhedral since $\gph\mathcal N_{D_i}$ is
	polyhedral for each $i\in L$ by \cref{lem:normal_cone_map_of_polyhedral_set}.
	Using invertibility of $A$, polyhedrality of $\gph\mathcal N_{D_i}$ for each $i\in L$,
	and \cref{lem:product_rule_tangents_polyhedral_sets}, however, we find
	\begin{align*}
		\mathcal T_{\gph\mathcal N_D}(\bar y,\bar y^*)
		&=
		\mathcal T_{A\prod_{i=1}^\ell\gph\mathcal N_{D_i}}(\bar y,\bar y^*)
		=
		A\mathcal T_{\prod_{i=1}^\ell\gph\mathcal N_{D_i}}((\bar y_1,\bar y_1^*),\ldots,(\bar y_\ell,\bar y_\ell^*))
		\\
		&=
		A\prod_{i=1}^\ell\mathcal T_{\gph\mathcal N_{D_i}}(\bar y_i,\bar y_i^*)
		\\
		&=
		\{(u,v)\in\R^m\times\R^m\,|\,\forall i\in L\colon\,(u_i,v_i)
		\in\mathcal T_{\gph\mathcal N_{D_i}}(\bar y_i,\bar y_i^*)\}
	\end{align*}
	which shows the formula for the graphical derivative.
\end{proof}

In the subsequently stated definition, we first recall the notion of the regular and limiting coderivative of a set-valued mapping
before introducing its so-called directional pseudo-coderivative. The latter will be of essential importance in the course of the paper.
Finally, we recall the notion of directional pseudo-coderivatives as introduced by Gfrerer in \cite[Definition~2]{Gfrerer2014a}.
\begin{definition}\label{def:coderivatives}
	Let $\Phi\colon\mathbb X\tto\mathbb Y$ be a set-valued mapping possessing a closed graph locally around $(\xb,\yb)\in\gph \Phi$.
	Furthermore, let $(u,v)\in\mathbb X\times\mathbb Y$ be a pair of directions.
	\begin{enumerate}
	\item The \emph{regular and limiting coderivative} of $\Phi$ at $(\bar x,\bar y)$
	are the set-valued mappings
	$\widehat D^*\Phi(\bar x,\bar y)\colon\mathbb Y\tto\mathbb X$ 
	and
	$D^*\Phi(\bar x,\bar y)\colon\mathbb Y\tto\mathbb X$ given,
	respectively, by
	\begin{align*}
		\forall y^*\in\mathbb Y\colon\quad
		\widehat D^*\Phi(\bar x,\bar y)(y^*)
		&:=
		\left\{x^*\in\mathbb X\,\middle|\,
			(x^*,-y^*)\in\widehat{\mathcal N}_{\gph\Phi}(\bar x,\bar y)\right\},
		\\
		D^*\Phi(\bar x,\bar y)(y^*)
		&:=
		\left\{x^*\in\mathbb X\,\middle|\,
			(x^*,-y^*)\in\mathcal N_{\gph\Phi}(\bar x,\bar y)\right\}.
	\end{align*}
	The set-valued mapping
	$D^*\Phi((\bar x,\bar y);(u,v))\colon\mathbb Y\tto\mathbb X$ given by
	\begin{align*}
		\forall y^*\in\mathbb Y\colon\quad
		D^*\Phi((\bar x,\bar y);(u,v))(y^*)
		&:=
		\left\{x^*\in\mathbb X\,\middle|\,
			(x^*,-y^*)\in\mathcal N_{\gph\Phi}((\bar x,\bar y);(u,v))\right\}
	\end{align*}
	is the \emph{limiting coderivative} of $\Phi$ at $(\bar x,\bar y)$ \emph{in direction} $(u,v)$.
	If $\Phi$ is single-valued at $\bar x$, we use 
	$\widehat D^*\Phi(\bar x),D^*\Phi(\bar x),D^*\Phi(\bar x;(u,v))\colon\mathbb Y\tto\mathbb X$
	for brevity.
	\item\label{item:new_pseudo_coderivative}
	Given $\gamma \geq 1$ and $u\in\mathbb S_{\mathbb X}$, the
 	\emph{pseudo-coderivative of order $\gamma$} of $\Phi$ at $(\bar x,\bar y)$ 
 	\emph{in direction $(u,v)$}
 	is the mapping $D^\ast_{\gamma} \Phi((\xb,\yb); (u,v))\colon\mathbb Y\tto\mathbb X$
 	which assigns to every $y^*\in\mathbb Y$ the set of all $x^*\in\mathbb X$ such that 
 	there are sequences $\{u_k\}_{k\in\N},\{x_k^*\}_{k\in\N}\subset\mathbb X$,
 	$\{v_k\}_{k\in\N},\{y_k^*\}_{k\in\N}\subset\mathbb Y$, 
 	and $\{t_k\}_{k\in\N}\subset\R_+$ which satisfy
 	$u_k\to u$, $v_k\to v$, $t_k\searrow 0$, $x_k^*\to x^*$, $y_k^*\to y^*$, and
 	\begin{align*}
 		\forall k\in\N\colon\quad
 		\left(x^\ast_k,-\frac{y^\ast_k}{(t_k\norm{u_k})^{\gamma-1}}\right)
		\in 
		\widehat{\mathcal N}_{\gph \Phi}(\xb + t_k u_k,\yb + (t_k\norm{u_k})^{\gamma} v_k).
 	\end{align*}
 	In case $\gamma:=1$, this definition recovers the one of $D^*\Phi((\bar x,\bar y);(u,v))$.
 	\item\label{item:Gfrerer_pseudo_coderivative} 
 	Given $\gamma\geq 1$ and $u\in\mathbb S_{\mathbb X}$, 
 	\emph{Gfrerer's pseudo-coderivative of order $\gamma$} of $\Phi$ 
 	at $(\bar x,\bar y)$ \emph{in direction $(u,v)$}
 	is the mapping 
 	$\widetilde D^\ast_{\gamma} \Phi((\xb,\yb); (u,v))\colon\mathbb Y\tto\mathbb X$
 	which assigns to every $y^*\in\mathbb Y$ the set of all $x^*\in\mathbb X$ 
 	such that there are sequences $\{u_k\}_{k\in\N},\{x_k^*\}_{k\in\N}\subset\mathbb X$,
 	$\{v_k\}_{k\in\N},\{y_k^*\}_{k\in\N}\subset\mathbb Y$, 
 	and $\{t_k\}_{k\in\N}\subset\R_+$ which satisfy
 	$u_k\to u$, $v_k\to v$, $t_k\searrow 0$, $x_k^*\to x^*$, $y_k^*\to y^*$, and
 	\begin{align*}
 		\forall k\in\N\colon\quad
 		\left(x^\ast_k,-\frac{y^\ast_k}{(t_k\norm{u_k})^{\gamma-1}}\right)
		\in 
		\widehat{\mathcal N}_{\gph \Phi}(\xb + t_k u_k,\yb + t_k v_k).
 	\end{align*}
 	Again, for $\gamma:=1$, we recover the definition of $D^*\Phi((\bar x,\bar y);(u,v))$.
	\end{enumerate}
\end{definition}

Let $\Phi\colon\mathbb X\tto\mathbb Y$ be a set-valued mapping whose graph is closed locally around $(\bar x,\bar y)\in\gph\Phi$
and fix a pair of directions $(u,v)\in\mathbb S_{\mathbb X}\times\mathbb Y$, $(x^*,y^*) \in\mathbb X\times\mathbb Y$, and $\gamma>1$.
Then we obtain the trivial relations
\begin{equation}\label{eq:trivial_upper_estimate_pseudo_coderivative}
	x^* \in D^*_\gamma\Phi((\bar x,\bar y);(u,v))(y^*)
	\ \Longrightarrow \
	\left\{\begin{aligned}
			&0 \in D\Phi(\bar x,\bar y)(u), \ 0 \in D^*\Phi(\xb,\yb)(y^*), \\
			&0 \in D^*\Phi((\xb,\yb);(u,0))(y^*),\\
			&x^* \in \widetilde D^*_\gamma\Phi((\bar x,\bar y);(u,0))(y^*).
		\end{aligned}
	\right.
\end{equation}

Graphical derivative and (directional) limiting coderivative 
are powerful tools for studying regularity properties of set-valued mappings,
such as (strong) metric regularity and subregularity, as well as their 
inverse counterparts of Lipschitzness,
such as Aubin property and (isolated) calmness.
Indeed, given a closed-graph set-valued mapping $\Phi\colon\mathbb X\tto\mathbb Y$, 
metric regularity and strong metric subregularity at some $(\xb,\yb) \in \gph \Phi$
are characterized, respectively, by
\[
	\ker D^*\Phi(\xb,\yb) = \{0\} \qquad \textrm{ and } \qquad \ker D\Phi(\xb,\yb) = \{0\},
\]
see e.g.\ \cite{Levy96,Mordukhovich2018,RockafellarWets1998} for the definition of
these Lipschitzian properties as well as the above results.
Moreover, metric subregularity of $\Phi$ at $(\xb,\yb)$ is implied by
\[
	\forall u\in\ker D\Phi(\xb,\yb)\cap\mathbb S_{\mathbb X}\colon\quad
	\ker D^*\Phi((\xb,\yb);(u,0)) = \{0\},
\]
which is referred to as \emph{First-Order Sufficient Condition for Metric Subregularity} 
(FOSCMS for short) in the literature,
see e.g.\ \cite{Gfrerer2013}.
The relations \eqref{eq:trivial_upper_estimate_pseudo_coderivative} suggest that
the pseudo-coderivative can be useful particularly in situations where the above regularity properties fail.

Note that the aforementioned notions of regularity and Lipschitzness 
express certain linear rate of change of the mapping.
Similarly, there is an underlying linearity in the definition 
of graphical derivative and coderivatives.
Take the graphical derivative for instance. 
Since the same sequence $\{t_k\}_{k\in\N}$ appears in the domain 
as well as in the range space, if $v \in D\Phi(\xb,\yb)(u)$ 
implies that $u\in\mathbb X$ and $v\in\mathbb Y$ are both nonzero,
it suggests a proportional (linear) rate of change.
Thus, in order to characterize pseudo (sub)regularity of order $\gamma>1$ of $\Phi$, 
it is not very surprising that 
we need to exploit derivative-like objects based on sub- or superlinear structure.
Exemplary, this has been successfully visualized in \cite[Corollary~2]{Gfrerer2014a} 
by means of Gfrerer's directional pseudo-coderivative
of order $\gamma>1$ from \cref{def:coderivatives}\,\ref{item:Gfrerer_pseudo_coderivative}.
Here, we show an analogous sufficient condition for metric pseudo-subregularity 
via the pseudo-coderivative from \cref{def:coderivatives}\,\ref{item:new_pseudo_coderivative}.

\begin{lemma}\label{lem:sufficient_condition_pseudo_subregularity_abstract}
	Let $\Phi\colon\mathbb X\tto\mathbb Y$ be a set-valued mapping 
	having locally closed graph around
	$(\bar x,\bar y)\in\gph\Phi$, 
	fix a direction $u\in\mathbb S_{\mathbb X}$, and some $\gamma\geq 1$.
	Assume that
	\begin{equation}\label{eq:FOSCMS_gamma}
		\ker D^\ast_{\gamma}  \Phi((\xb,\yb ); (u,0))=\{0\}
	\end{equation}
	holds.
	Then $\Phi$ is metrically pseudo-subregular of order $\gamma$ 
	at $(\bar x,\bar y)$ in direction $u$.
\end{lemma}
\begin{proof}
	Suppose that $\Phi$ is not metrically pseudo-subregular 
	of order $\gamma$ at $(\bar x,\bar y)$
	in direction $u\in\mathbb Y$. 
	Due to \cite[Theorem~1(2)]{Gfrerer2014a}, we find sequences $\{t_k\}_{k\in\N}\subset\R_+$,
	$\{u_k\}_{k\in\N},\{x_k^*\}_{k\in\N}\subset\mathbb X$, 
	and $\{v_k\}_{k\in\N},\{y_k^*\}_{k\in\N}\subset\mathbb Y$
	satisfying (among other things) $t_k\searrow 0$, $u_k\to u$, $t_k^{1-\gamma}v_k\to 0$, 
	$x_k^*\to 0$, such that
	$\norm{y_k^*}=1$ and 
	\[
		(x_k^*,-y_k^*/(t_k\norm{u_k})^{\gamma-1})
		\in
		\widehat{\mathcal N}_{\gph\Phi}((\bar x,\bar y)+t_k(u_k,v_k))
	\]
	for each $k\in\N$. Let us set $\tilde v_k:=t_k^{1-\gamma}\norm{u_k}^{-\gamma} v_k$ 
	for each $k\in\N$.
	Then we have 
	\[
		(x_k^*,-y_k^*/(t_k\norm{u_k})^{\gamma-1})
		\in
		\widehat{\mathcal N}_{\gph\Phi}(\bar x+t_ku_k,\bar y+(t_k\norm{u_k})^\gamma \tilde v_k)
	\]
	for each $k\in\N$ and $\tilde v_k\to 0$ from $t_k^{1-\gamma}v_k\to 0$.
	Observing that $\{y_k^*\}_{k\in\N}$ possesses 
	a nonvanishing accumulation point $y^*\in\mathbb Y$,
	taking the limit along a suitable subsequence 
	yields $0\in  D^*_\gamma((\bar x,\bar y);(u,0))(y^*)$
	which contradicts the assumptions of the lemma.
\end{proof}

Let us remark that due to \eqref{eq:trivial_upper_estimate_pseudo_coderivative}, condition
\begin{equation}\label{eq:FO_characterization_of_dir_metric_pseudo_reg}
	\ker \widetilde{D}^*_\gamma\Phi((\bar x,\bar y);(u,0))=\{0\}
\end{equation}
is stronger than \eqref{eq:FOSCMS_gamma} and, thus, also sufficient 
for metric pseudo-subregularity
of $\Phi$ of order $\gamma\geq 1$ at $(\bar x,\bar y)$ in direction $u$.
By means of \cite[Corollary~2]{Gfrerer2014a},
\eqref{eq:FO_characterization_of_dir_metric_pseudo_reg}
is actually equivalent to $\Phi$ being so-called 
metrically pseudo-regular at $(\bar x,\bar y)$ in
direction $(u,0)$, see \cite[Definition~1]{Gfrerer2014a} for a definition.
Note that in case $\gamma:=1$, both conditions \eqref{eq:FOSCMS_gamma} 
and \eqref{eq:FO_characterization_of_dir_metric_pseudo_reg}
recover FOSCMS of $\Phi$ at $(\bar x,\bar y)$ in direction $u$.

\section{Variational analysis of constraint mappings}\label{sec:variational_analysis_constraint_mapping}

In this section, we address the pseudo-coderivative calculus for constraint mappings.
Therefore, let us fix a twice continuously differentiable mapping 
$g\colon\mathbb X\to\mathbb Y$ as well as
a closed set $D\subset\mathbb Y$ and consider the mapping $\Phi\colon\mathbb X\tto\mathbb Y$
given by $\Phi(x):=g(x)-D$ for all $x\in\mathbb X$.
To start, let us consider the regular and (directional) limiting coderivative of $\Phi$
which can be easily computed by applying the change-of-coordinates formulas provided in
\cite[Exercise~6.7]{RockafellarWets1998} and \cite[Lemma~2.1]{BenkoMehlitz2020}.

\begin{lemma}\label{lem:coderivatives_constraint_maps}
	Fix $(x,y)\in\gph\Phi$. Then, for each $y^*\in\mathbb Y$
	and each pair $(u,v)\in\mathbb X\times\mathbb Y$, we have
    \begin{align*}
		\widehat{D}^*\Phi(x,y)(y^*)
		&=
		\begin{cases}
			\nabla g(x)^*y^*	&	y^*\in \widehat{\mathcal N}_D(g(x)-y),\\
			\varnothing				&	\text{otherwise},
		\end{cases}
		\\
		D^*\Phi(x,y)(y^*)
		&=
		\begin{cases}
			\nabla g(x)^*y^*	&	y^*\in \mathcal N_D(g(x)-y),\\
			\varnothing	&	\text{otherwise},
		\end{cases}
		\\
		D^*\Phi((x,y);(u,v))(y^*)
		&=
		\begin{cases}
			\nabla g(x)^*y^*	& y^*\in\mathcal N_D(g(x)-y;\nabla g(x)u-v),\\
			\varnothing				&\text{otherwise.}
		\end{cases}
	\end{align*}
\end{lemma}

Next, we estimate the directional pseudo-coderivatives of order $2$ 
of constraint mappings in terms of initial problem data.
The technical proof of the subsequently stated result can be
found in \cref{sec:appendix}.

\begin{theorem}\label{The : NCgen}
Given $(\xb,0) \in \gph \Phi$ and a direction $u \in \mathbb S_{\mathbb X}$, let
\[
	x^* \in \widetilde{D}^\ast_{2} \Phi((\xb,0);(u,v))(y^*)
\]
for some $v, y^*\in\mathbb Y$.
Then the following statements hold.
\begin{enumerate}
\item\label{item:general_estimate}
There exists $z^*\in[y^*]^\perp$ such that
\begin{subequations}\label{eq:upper_estimate_pseudo_coderivative_order_two}
	\begin{align}
		\label{eq:2ordEstimNC}
 		x^* &= \nabla^2\langle y^*,g\rangle(\bar x)(u) + \nabla g(\xb)^* z^*,
 		\\
 		\label{eq:multiplier_from_dir_lim_normal_cone_and_kernel}
		y^* &\in \mathcal N_{D}(g(\xb);\nabla g(\xb) u - v) \cap \ker \nabla g(\xb)^*.
	\end{align}
\end{subequations}
\item\label{item:general_estimate_+CQ}
Each of the following two conditions
\begin{subequations}\label{eq:some:CQ}
	\begin{align}
		\label{eq:some:CQ_1}
 		& D\mathcal N_{D}(g(\xb),y^*)(0) \cap \ker \nabla g(\xb)^*=\{0\},
 		\\
 		\label{eq:some:CQ_2}
		\nabla g(\xb)u \neq v,
		\quad
		& D_{\textup{sub}}\mathcal N_{D}(g(\xb),y^*)
		\left(\frac{\nabla g(\xb)u - v}{\norm{\nabla g(\xb)u - v}} \right) 
		\cap \ker \nabla g(\xb)^*
		=\varnothing,
	\end{align}
\end{subequations}
implies that we can find $z^* \in D\mathcal N_{D}(g(\xb),y^*)(\nabla g(\xb)u - v)$ satisfying \eqref{eq:upper_estimate_pseudo_coderivative_order_two}.
\item\label{item:polyhedral_estimate}
If $\mathbb Y:=\R^m$ and $D$ is locally polyhedral around $g(\bar x)$, 
then $\mathcal N_{D}(g(\xb);\nabla g(\xb) u - v)= \mathcal N_{\mathcal T_D(g(\xb))}(\nabla g(\xb) u - v)$ 
and there exists
$z^*\in\mathcal T_{\mathcal N_{\mathcal T_D(g(\xb))}(\nabla g(\xb) u - v)}(y^*)$
satisfying \eqref{eq:upper_estimate_pseudo_coderivative_order_two}.
\item\label{item:polyhedral_estimate_order_2}
If $\mathbb Y:=\R^m$, $D$ is locally polyhedral around $g(\bar x)$, and we even have
\[
	x^* \in D^\ast_{2} \Phi((\xb,0);(u,v))(y^*) 
\]
for some $v, y^*\in\R^m$, there exist $s\in\mathbb X$ and 
$z^*\in\mathcal T_{\mathcal N_{\mathbf T(u)}(w_s(u,v))}(y^*)$
satisfying \eqref{eq:2ordEstimNC} and 
$y^* \in \mathcal N_{\mathbf T(u)}(w_s(u,v))\cap \ker \nabla g(\bar x)^*$ where
\begin{equation}\label{eq:Tu_and_ws}
 \mathbf T(u):=\mathcal T_{\mathcal T_D(g(\xb))}(\nabla g(\xb) u),\qquad
 w_s(u,v):= \nabla g(\xb) s + 1/2 \nabla^2 g(\xb)[u,u] - v.
\end{equation}
\end{enumerate}
\end{theorem}

In the following remark, we comment on the findings of \cref{The : NCgen}.
\begin{remark}\label{rem:upper_estimate_pseudo_coderivative_constraint_maps}
	We use the notation from \cref{The : NCgen}.
	\begin{enumerate}
	\item\label{item:comparison_estimate_pseudocoderivative_polyhedral}
	Note the similarity between cases~\ref{item:general_estimate_+CQ} and \ref{item:polyhedral_estimate} of \cref{The : NCgen}.
	If $D$ is actually a \emph{convex} polyhedral set in $\mathbb Y:=\R^m$, 
	based on the so-called reduction lemma, see \cite[Lemma~2E.4]{DontchevRockafellar2014}, 
	for each pair $(\bar z, \bar z^*) \in \gph \mathcal N_D$, we get
	\begin{align*}
		&\big(\gph \mathcal N_D - (\bar z,\bar z^*)\big) \cap \mathcal{O} 
		\\
		&\qquad
		=
		 \{(w,w^*) \,|\, w \in \mathcal K_D(\bar z,\bar z^*),\, w^* \in \mathcal K_D(\bar z,\bar z^*)^{\circ},\, \innerprod{w}{ w^* }=0\} \cap \mathcal{O},
	\end{align*}
	where $\mathcal{O}\subset\R^m\times\R^m$ is a neighborhood of $0$ and 
	$\mathcal K_D(\bar z,\bar z^*) := \mathcal T_D(\bar z) \cap [\bar z^*]^\perp$ stands for the \emph{critical cone} to $D$ at $(\bar z,\bar z^*)$,
	see \cite[Lemma~3.8.]{Benko2019} as well. 
	By \cref{lem:some_properties_of_polyhedral_sets}\,\ref{item:exactness_tangential_approximation}, this simply means
	\[
		\mathcal T_{\gph \mathcal N_D}(\bar z,\bar z^*) 
		=
		\{(w,w^*) \,|\, w \in \mathcal K_D(\bar z,\bar z^*),\, w^* \in \mathcal K_D(\bar z,\bar z^*)^{\circ},\, \innerprod{w}{ w^* }=0\}.
	\]
	Thus, $z^* \in D\mathcal N_{D}(g(\xb),y^*)(\nabla g(\xb)u - v)$ 
	means $\nabla g(\xb)u - v\in \mathcal T_D(g(\xb)) \cap [y^*]^\perp$, which gives us
	\begin{align*}
		y^*
		\in \mathcal N_D(g(\bar x))\cap[\nabla g(\bar x)u-v]^\perp
		= 
		\mathcal N_{\mathcal T_D(g(\xb))}(\nabla g(\xb) u - v),
	\end{align*}
	and
	\begin{align*}
		z^*
		& \in 
		\big(\mathcal T_D(g(\xb)) \cap [y^*]^\perp\big)^{\circ} \cap [\nabla g(\xb)u - v]^\perp\\
		& =
		\big(\mathcal N_D(g(\xb))+\spa(y^*)\bigr)\cap [\nabla g(\bar x)u-v]^\perp\\
		& =
		\mathcal N_D(g(\bar x))\cap[\nabla g(\bar x)u-v]^\perp + \spa(y^*)\\
		& = 
		\mathcal N_{\mathcal T_D(g(\xb))}(\nabla g(\xb)u - v) + \spa(y^*)
		= 
		\mathcal T_{\mathcal N_{\mathcal T_D(g(\xb))}(\nabla g(\xb) u - v)}(y^*)
	\end{align*}
	by the basic properties of convex polyhedral cones and \cref{lem:some_properties_of_polyhedral_sets}\,\ref{item:normal_cones_to_polyhedral_sets}.
	Hence, in the convex polyhedral case, the information on $y^*$ and $z^*$ from statements~\ref{item:general_estimate_+CQ} and \ref{item:polyhedral_estimate}
	of \cref{The : NCgen} is the same.
	\item
	The assumptions \eqref{eq:some:CQ} in \cref{The : NCgen}\,\ref{item:general_estimate_+CQ} rule out the possibility 
	that the sequence $\{(y_k^* - y^*)/\tau_k\}_{k \in \N}$ from the proof (see \cref{sec:appendix})
	is unbounded, while, in \cref{The : NCgen}\,\ref{item:polyhedral_estimate}, the underlying polyhedrality just ensures the existence of a bounded sequence,
	which can be used instead of $\{(y_k^* - y^*)/\tau_k\}_{k \in \N}$.
	\item 
	Note that, in case $\nabla g(\bar x)u\neq v$, assumption \eqref{eq:some:CQ_2}, 
	which is stated in terms of the graphical subderivative,
	is milder than \eqref{eq:some:CQ_1} in terms or the standard graphical derivative,
	and it preserves the connection to the direction $\nabla g(\xb)u - v$.
	Let us also note that the case $\nabla g(\bar x)u=v$ is, anyhow, special since this
	would annihilate the directional information in \eqref{eq:multiplier_from_dir_lim_normal_cone_and_kernel}
	completely.
	\item\label{item:graphical_subderivative_assumption_polyhedral}
	If $\mathbb Y:=\R^m$ and $D$ is locally polyhedral around $g(\bar x)$,
	conditions \eqref{eq:some:CQ} reduce to
	\begin{equation}\label{eq:CQ_pseudo_subregularity_I_polyhedral}
		D\mathcal N_{D}(g(\xb),y^*)(\nabla g(\xb)u - v) 
		\cap \ker \nabla g(\xb)^*
		\subset \{0\}
	\end{equation}
	thanks to \cref{lem:normal_cone_map_of_polyhedral_set}.
	\item 
	Let us mention that if $\mathbb Y:=\R^m$ and $D$ is locally polyhedral around $g(\bar x)$,
	we get the relations 
	\begin{align*}
            \mathcal N_{\mathbf T(u)}(w_s(u,v))
            & = 
            \mathcal N_{\mathcal T_D(g(\bar x))}(\nabla g(\bar x)u;w_s(u,v))
            \\
            & \subset 
            \mathcal N_{\mathcal T_D(g(\bar x))}(\nabla g(\bar x)u)\cap[w_s(u,v)]^\perp\\
            & \subset 
            \mathcal N_{\mathcal T_D(g(\bar x))}(\nabla g(\bar x)u)
             = 
            \mathcal N_D(g(\bar x);\nabla g(\bar x)u).
     \end{align*}
     from 
     \cref{lem:some_properties_of_polyhedral_sets}\,\ref{item:normal_cones_to_polyhedral_sets}.   
     This directly yields 
     $\mathcal T_{\mathcal N_{\mathbf T(u)}(w_s(u,v))}(y^*) 
     \subset\mathcal T_{\mathcal N_{\mathcal T_D(g(\bar x))}(\nabla g(\bar x)u)}(y^*)$.
	\end{enumerate}
\end{remark}

Taking into account the sufficient conditions for metric pseudo-subregularity from \cref{lem:sufficient_condition_pseudo_subregularity_abstract},
we obtain the following result for constraint mappings.

\begin{corollary}\label{cor:sufficient_condition_pseudo_subregularity}
 Given $(\xb,0) \in \gph \Phi$ and a direction $u \in \mathbb S_{\mathbb X}$,
 $\Phi$ is metrically pseudo-subregular of order 2 in direction $u$
 at $(\xb,0)$, provided one of the following conditions holds.
 \begin{enumerate}
 	 \item One has 
  		\begin{equation*}
   			\left. \begin{aligned}
        		&\nabla g(\xb)^* y^* = 0, \, 
        		\nabla^2\langle y^*,g\rangle(\bar x)(u) + \nabla g(\xb)^* z^* = 0, \\
        		&y^* \in \mathcal N_{D}(g(\xb);\nabla g(\xb) u), \,
        		 z^* \in [y^*]^{\perp} 
       			\end{aligned}
				\right\} 
				\quad \Longrightarrow \quad 
				y^* = 0.
  		\end{equation*}
  	\item One has
  		\begin{equation}\label{eq:CQ_pseudo_subregularity_II}
  			\left. \begin{aligned}
        		&\nabla g(\xb)^* y^* = 0, \, 
        		\nabla^2\langle y^*,g\rangle(\bar x)(u) + \nabla g(\xb)^* z^* = 0, \\
        		&y^* \in \mathcal N_{D}(g(\xb);\nabla g(\xb) u), \,
        		z^*\in D\mathcal N_D(g(\bar x),y^*)(\nabla g(\bar x)u)
       		\end{aligned}
       		\right\} 
       		\quad \Longrightarrow \quad 
       		y^* = 0.
  		\end{equation}
  		Furthermore, we either have
  		\begin{equation}\label{eq:CQ_pseudo_subregularity_Ia}
  			\left. \begin{aligned}
        		&\nabla g(\xb)^* y^* = 0, \, 
        		\nabla g(\xb)^* \hat z^* = 0, \\
        		&y^* \in \mathcal N_{D}(g(\xb);\nabla g(\xb) u), \,
        		 \hat z^* \in D\mathcal N_D(g(\bar x),y^*)(0) 
       			\end{aligned}
				\right\} 
				\quad \Longrightarrow \quad 
				\hat z^* = 0
  		\end{equation}
		or $\nabla g(\bar x)u\neq 0$ and
		\begin{equation}\label{eq:CQ_pseudo_subregularity_Ib}
  			\left. \begin{aligned}
        		&\nabla g(\xb)^* y^* = 0, \, 
        		\nabla g(\xb)^* \hat z^* = 0, \\
        		&y^* \in \mathcal N_{D}(g(\xb);\nabla g(\xb) u)
       			\end{aligned}
				\right\} 
				\ \Longrightarrow \ 
				\hat z^* \notin D_\textup{sub}\mathcal N_D(g(\bar x),y^*)
				\left(\frac{\nabla g(\bar x)u}{\norm{\nabla g(\bar x)u}}\right).
  		\end{equation}
  		\item It holds $\mathbb Y:=\R^m$, $D$ is locally polyhedral around $g(\bar x)$,
  			and one has
  			\begin{equation}\label{eq:CQ_pseudo_subregularity_polyhedral_I}
  			\left. \begin{aligned}
        		&\nabla g(\xb)^* y^* = 0, \, 
        		\nabla^2\langle y^*,g\rangle(\bar x)(u) + \nabla g(\xb)^* z^* = 0, \\
        		&y^* \in \mathcal N_{\mathcal T_D(g(\bar x))}(\nabla g(\xb) u), \,
        		z^*\in \mathcal T_{\mathcal N_{\mathcal T_D(g(\bar x))}(\nabla g(\bar x)u)}(y^*)
       		\end{aligned}
       		\right\} 
       		\quad \Longrightarrow \quad 
       		y^* = 0.
  			\end{equation}
 		\item It holds $\mathbb Y:=\R^m$, $D$ is locally polyhedral around $g(\xb)$, 
 			and one has
 			\begin{equation}\label{eq:CQ_pseudo_subregularity_polyhedral_II}
   				\left. \begin{aligned}
        			&\nabla g(\xb)^* y^* = 0, \, 
        			\nabla^2\langle y^*,g\rangle(\bar x)(u) + \nabla g(\xb)^* z^* = 0, \\
        			&y^* \in \mathcal N_{\mathbf T(u)}(w_s(u,0)), \, 
        			z^* \in \mathcal T_{\mathcal N_{\mathbf T(u)}(w_s(u,0))}(y^*)  
       				\end{aligned}
				\right\} 
				\quad \Longrightarrow \quad y^* = 0
  			\end{equation}
  			for each $s\in\mathbb X$,
  			where $\mathbf T(u)$ and $w_s(u,0)$ are given as in \eqref{eq:Tu_and_ws}.
 		\end{enumerate}
\end{corollary}

We would like to point the reader's attention to \cite[Theorem 2]{Gfrerer2014a}
where related results were obtained. Let us note, however, that the conditions
stated in the first
and second assertion of \cref{cor:sufficient_condition_pseudo_subregularity} are not
covered by \cite[Theorem~2]{Gfrerer2014a} since $D$ does not need
to be polyhedral for our findings.
We note that \eqref{eq:CQ_pseudo_subregularity_polyhedral_I} is a slightly stronger
condition than the one derived in \cite[Theorem~2]{Gfrerer2014a} where 
an additional extremal condition w.r.t.\ the active components
of the polyhedral set $D$ is incorporated.
Without going into details, it seems that this extremal condition
can be similarly strong as the requirement $y^* \in \mathcal N_{\mathbf T(u)}(w_s(u,0))$,
since the duality theory of linear programming appears to yield $\langle y^*,w_s(u,0)\rangle = 0$.
This would put the conditions from \cite[Theorem~2]{Gfrerer2014a} roughly somewhere between
our conditions \eqref{eq:CQ_pseudo_subregularity_polyhedral_I} and \eqref{eq:CQ_pseudo_subregularity_polyhedral_II}.

For $D:=\{0\}$, conditions \eqref{eq:CQ_pseudo_subregularity_polyhedral_I} and \eqref{eq:CQ_pseudo_subregularity_polyhedral_II}
coincide with
\[
	\nabla g(\bar x)^*y^*=0,\,\nabla^2\langle y^*,g\rangle(\bar x)(u)+\nabla g(\bar x)^*z^*=0
	\quad\Longrightarrow\quad
	y^*=0.
\]
The latter is equivalent to so-called $2$-regularity of the function $g$ at $\bar x$ in direction $u$,
see \cite[Proposition~2]{GfrererOutrata2016c} for this result and 
\cite{Avakov1985} for the origins of this property.
It can also be distilled from \cite[Proposition~2]{GfrererOutrata2016c} that $2$-regularity
of $g$ at $\bar x$ in direction $u$ is equivalent to metric pseudo-regularity of $g$ at $(\bar x,g(\bar x))$
in direction $(u,\nabla g(\bar x)u)$ of order $2$.
Furthermore, we want to point the reader's attention to the fact that there is a connection between
\eqref{eq:CQ_pseudo_subregularity_polyhedral_I} or \eqref{eq:CQ_pseudo_subregularity_polyhedral_II}
and $2$-regularity of $g$ at $\bar x$ w.r.t.\ $D$ in direction $u$
as introduced in \cite[Definition~1]{ArutyunovAvakovIzmailov2008}, and the latter property is, thus,
related to directional pseudo-(sub)regularity of constraint mappings of order $2$.
A precise comparison of those properties is, however, beyond the capacity of this paper.

\section{Stationarity conditions for nonsmooth optimization problems involving pseudo-coderivatives}\label{sec:directional_asymptotic_tools}

\subsection{General case}

For a locally Lipschitz continuous function $\varphi\colon\mathbb X\to\R$
and a set-valued mapping $\Phi\colon\mathbb X\tto\mathbb Y$ with a closed graph
and $\bar y\in\Im\Phi$, we investigate the rather abstract optimization problem
\begin{equation}\label{eq:nonsmooth_problem}\tag{P}
	\min\{\varphi(x)\,|\,\bar y\in\Phi(x)\}.
\end{equation}
Throughout the section, the feasible set of \eqref{eq:nonsmooth_problem} 
will be denoted by $\mathcal F\subset\mathbb X$. 
Clearly, we have $\mathcal F\neq\varnothing$ from $\bar y\in\Im\Phi$.
Note that the model \eqref{eq:nonsmooth_problem} covers numerous 
classes of optimization problems from the literature including
standard nonlinear problems, problems with geometric 
(particularly, disjunctive or conic) constraints, problems with
(quasi-) variational inequality constraints, and bilevel optimization problems. 
Furthermore, we would like to mention that choosing $\bar y:=0$ would not be restrictive
since one could simply consider $\widetilde\Phi\colon\mathbb X\tto\mathbb Y$ given by
$\widetilde{\Phi}(x):=\Phi(x)-\bar y$, $x\in\mathbb X$, in case where $\bar y$ does not
vanish.
A standard notion of stationarity which applies to \eqref{eq:nonsmooth_problem} is the
one of M-stationarity.

\begin{definition}\label{def:M_stationarity}
	A feasible point $\bar x\in\mathcal F$ is called \emph{M-stationary} whenever there is a multiplier $\lambda\in\mathbb Y$ such that
	\[
		0\in\partial\varphi(\bar x)+D^*\Phi(\bar x,\bar y)(\lambda).
	\]
\end{definition}

Let us also recall a quite general notion of critical directions associated with \eqref{eq:nonsmooth_problem} taken from
\cite[Definition~5]{Gfrerer2013}.
\begin{definition}\label{def:critical_direction}
	For some feasible point $\bar x\in\mathcal F$, 
	a direction $u\in\mathbb S_{\mathbb X}$ is called \emph{critical}
	for \eqref{eq:nonsmooth_problem} at $\bar x$ whenever there are sequences $\{u_k\}_{k\in\N}\subset\mathbb X$, $\{v_k\}_{k\in\N}\subset\mathbb Y$,
	and $\{t_k\}_{k\in\N}\subset\R_+$ such that $u_k\to u$, $v_k\to 0$, $t_k\searrow 0$, and $(\bar x+t_ku_k,\bar y+t_kv_k)\in\gph\Phi$ for all $k\in\N$ as well as
	\[
		\limsup\limits_{k\to\infty}\frac{\varphi(\bar x+t_ku_k)-\varphi(\bar x)}{t_k}\leq 0.
	\]
\end{definition}

Clearly, each critical direction $u\in\mathbb S_{\mathbb X}$ of \eqref{eq:nonsmooth_problem} at $\bar x\in\mathcal F$ 
satisfies $u\in\ker D\Phi(\bar x,\bar y)$ by definition of the graphical derivative.
In case where $\varphi$ is, additionally, assumed to be directionally differentiable 
at $\bar x$, it is easily seen that
$u\in\mathbb S_{\mathbb X}$ is critical for \eqref{eq:nonsmooth_problem} at $\bar x$ if and only if
$\varphi'(\bar x;u)\leq 0$ and $0\in D\Phi(\bar x,\bar y)(u)$,
see \cite[Proposition~3.5]{Shapiro1990} as well.
Let us note that whenever $\bar x\in\mathcal F$ is a feasible point of \eqref{eq:nonsmooth_problem}
such that no critical direction for \eqref{eq:nonsmooth_problem} at $\bar x$ exists, then $\bar x$
is a strict local minimizer of \eqref{eq:nonsmooth_problem}.
Conversely, there may exist strict local minimizers of \eqref{eq:nonsmooth_problem} such that
a critical direction for \eqref{eq:nonsmooth_problem} at this point exists.

The following result first provides necessary optimality conditions for \eqref{eq:nonsmooth_problem}
which hold in the absence of constraint qualifications and are stated with the aid of the
pseudo-coderivative construction from \cref{def:coderivatives}.
More precisely, the derived conditions depend on a certain \emph{order} $\gamma\geq 1$.
Furthermore, the theorem specifies how metric pseudo-subregularity of $\Phi$ of order $\gamma$
at the reference point enriches the situation.

\begin{theorem}\label{thm:higher_order_directional_asymptotic_stationarity}
	Let $\bar x\in\mathcal F$ be a local minimizer of \eqref{eq:nonsmooth_problem} 
	and consider $\gamma \geq 1$.
	Then one of the following conditions holds.
	\begin{enumerate}
	  \item\label{item:trivial_MSt}
	  	The point $\bar x$ is M-stationary for \eqref{eq:nonsmooth_problem}.
	  \item\label{item:pseudo_MSt_zero_mult}%
		There exists a critical direction $u\in\mathbb S_{\mathbb X}$ 
		for \eqref{eq:nonsmooth_problem} at $\bar x$ such that
			\begin{equation}\label{eq:Pseudo-M-Stat_Gferer}
				0\in \partial\varphi(\bar x)+\widetilde{D}^*_\gamma\Phi((\bar x,\bar y);(u,0))(0).
			\end{equation}
	  \item\label{item:pseudo_dir_MSt}%
		There exist a critical direction $u\in\mathbb S_{\mathbb X}$ 
		for \eqref{eq:nonsmooth_problem} at $\bar x$
		and $v \in \mathbb Y$ as well as a nonvanishing multiplier $\lambda \in \mathbb Y$ such that
		$\innerprod{\lambda}{ v }\geq 0$ and
		\begin{equation}\label{eq:Pseudo-M-stat}	
		  0 \in \partial\varphi(\bar x) + D^\ast_{\gamma} \Phi((\xb,\yb);(u,v))(\lambda).
		\end{equation}
	  \item\label{item:irregular_case_pseudo_MSt}
		There exist a critical direction $u\in\mathbb S_{\mathbb X}$ 
		for \eqref{eq:nonsmooth_problem} at $\bar x$ 
	and $\{x_k\}_{k\in\N},\{x_k'\}_{k\in\N},\{\eta_k\}_{k\in\N}\subset\mathbb X$ as well as
	$\{y_k\}_{k\in\N},\{y_k^*\}_{k\in\N}\subset\mathbb Y$ such that
	$x_k,x_k'\notin\Phi^{-1}(\bar y)$, $y_k\neq\bar y$, and $y_k^*\neq 0$ for all $k\in\N$,
	satisfying the convergence properties 
	\begin{equation}\label{eq:convergences_directional_asymptotic_stationarity}
		\begin{aligned}
			x_k&\to\bar x,&		\qquad	x_k'&\to\bar x,&	\qquad	y_k&\to\bar y,&
			\\
			\frac{x_k-\bar x}{\nnorm{x_k-\bar x}}&\to u,&	\qquad
			\frac{x_k'-\bar x}{\nnorm{x_k'-\bar x}}&\to u,&	\qquad
			v_k^\gamma&\to 0,&
			\\
			y_k^*&\to y^*,&\qquad
			\nnorm{\lambda_k^\gamma} &\to \infty,&\qquad
			\eta_k&\to 0,&
			\\
			\frac{y_k-\bar y}{\nnorm{y_k-\bar y}}-\frac{\lambda_k^\gamma}{\nnorm{\lambda_k^\gamma}}&\to 0,& && 	&&
		\end{aligned}
	\end{equation}
	where we used
	\begin{equation}\label{eq:definition_of_surrogate_sequences}
		\forall k\in\N\colon\quad
		v^\gamma_k := \frac{y_k-\bar y}{\norm{x_k-\bar x}^\gamma},\qquad
		\lambda_k^\gamma:= \frac{\norm{x_k - \bar x}^\gamma}{\norm{y_k - \bar y}}y_k^*,
	\end{equation}
	as well as
	\begin{equation}\label{eq:asymptotic_stationarity_regular_tools_gamma}
		\forall k\in\N\colon\quad
		\eta_k\in\widehat{\partial}\varphi(x_k')+\widehat{D}^*\Phi(x_k,y_k)\left(\frac{\lambda_k^\gamma}{\norm{x_k - \bar x}^{\gamma-1}}\right).
	\end{equation}
	\end{enumerate}
	Moreover, if $\Phi$ is metrically pseudo-subregular of order $\gamma$ at $(\bar x, \bar y)$
	in each direction $u\in\mathbb S_{\mathbb X}\cap\ker D\Phi(\bar x,\bar y)$, 
	$\bar x$ satisfies one of the alternatives~\ref{item:trivial_MSt},~\ref{item:pseudo_MSt_zero_mult}, or~\ref{item:pseudo_dir_MSt}.
\end{theorem}
\begin{proof}
	Let $\varepsilon>0$ be chosen such that $\varphi(x)\geq\varphi(\bar x)$ holds for all $x\in\mathcal F\cap\mathbb B_\varepsilon(\bar x)$
	and, for each $k\in\N$, consider the optimization problem
	\begin{equation}\label{eq:penalized_nonsmooth_problem}\tag{P$(k)$}
		\min\limits_{x,y}\{\varphi(x)+\tfrac k2\norm{y-\bar y}^2+\tfrac12\norm{x-\bar x}^2\,|\,(x,y)\in\gph\Phi,\,x\in\mathbb B_\varepsilon(\bar x)\}.
	\end{equation}
	For each $k\in\N$, the objective function of \eqref{eq:penalized_nonsmooth_problem} is bounded from below, continuous on the closed
	feasible set of this problem, and coercive in the variable $y$,
	so \eqref{eq:penalized_nonsmooth_problem} possesses a global minimizer $(\bar x_k,\bar y_k)\in\mathbb X\times\mathbb Y$. 
	By feasibility of $(\bar x,\bar y)$ for
	\eqref{eq:penalized_nonsmooth_problem}, we find
	\begin{equation}\label{eq:estimate_from_penalization}
		\forall k\in\N\colon\quad 
		\varphi(\bar x_k)+\tfrac{k}{2}\norm{\bar y_k-\bar y}^2+\tfrac12\norm{\bar x_k-\bar x}^2\leq \varphi(\bar x).
	\end{equation}

	By boundedness of $\{\bar x_k\}_{k\in\N}\subset\mathbb B_\varepsilon(\bar x)$,
	we may assume $\bar x_k\to\tilde x$ for some $\tilde x\in\mathbb B_\varepsilon(\bar x)$.
	Observing that $\{\varphi(\bar x_k)\}_{k\in\N}$ is bounded by continuity of $\varphi$, $\bar y_k\to\bar y$ easily follows from
	\eqref{eq:estimate_from_penalization}. Furthermore, the closedness of $\gph\Phi$ guarantees $(\tilde x,\bar y)\in\gph\Phi$, i.e.,
	$\tilde x\in\mathcal F\cap\mathbb B_\varepsilon(\bar x)$ leading to $\varphi(\bar x)\leq\varphi(\tilde x)$.
	From \eqref{eq:estimate_from_penalization}, we find
	\begin{align*}
		\varphi(\bar x)
		\leq
		\varphi(\tilde x)
		\leq
		\varphi(\tilde x)+\tfrac12\norm{\tilde x-\bar x}^2
		=
		\lim\limits_{k\to\infty}\left(\varphi(\bar x_k)+\tfrac12\norm{\bar x_k-\bar x}^2\right)
		\leq
		\varphi(\bar x),
	\end{align*}
	and $\tilde x=\bar x$ follows. Thus, we have $\bar x_k\to\bar x$.
	
	Let us assume that there is some $k_0\in\N$ such that $\bar x_{k_0}$ is feasible to
	\eqref{eq:nonsmooth_problem}. 
	By \eqref{eq:estimate_from_penalization}, we find 
	\begin{align*}
		&\varphi(\bar x) +\tfrac{k_0}{2}\norm{\bar y_{k_0}-\bar y}^2+\tfrac12\norm{\bar x_{k_0}-\bar x}^2
		\\
		&\qquad
		\leq 
		\varphi(\bar x_{k_0})+\tfrac{k_0}{2}\norm{\bar y_{k_0}-\bar y}^2+\tfrac12\norm{\bar x_{k_0}-\bar x}^2
		\leq
		\varphi(\bar x),
	\end{align*}
	i.e., $\bar x_{k_0}=\bar x$ and $\bar y_{k_0}=\bar y$. 
	Applying \cite[Theorem~6.1]{Mordukhovich2018}, 
	the subdifferential sum rule \cite[Theorem~2.19]{Mordukhovich2018}, and the
	definition of the limiting coderivative to find stationarity conditions
	of \eqref{eq:penalized_nonsmooth_problem} at $(\bar x,\bar y)$ yields
	$0\in\partial\varphi(\bar x)+D^*\Phi(\bar x,\bar y)(0)$,
	which is covered by~\ref{item:trivial_MSt}.

	Thus, we may assume that $\bar x_k\notin\mathcal F=\Phi^{-1}(\bar y)$ holds for all $k\in\N$. 
	Particularly, $\bar x_k\neq\bar x$ and $\bar y_k\neq\bar y$ is valid for all $k\in\N$ in this situation.	
	Assume without loss of generality that $\{\bar x_k\}_{k\in\N}$ belongs to the interior
	of $\mathbb B_\varepsilon(\bar x)$.

	We can apply Fermat's rule, see \cite[Proposition~1.30~(i)]{Mordukhovich2018}, 
	the fuzzy sum rule for regular normals from \cite[Exercise~2.26]{Mordukhovich2018}, and the
	definition of the regular coderivative in order to find sequences
	$\{(x_k,y_k)\}_{k\in\N}\subset\gph\Phi$, $\{x_k'\}_{k\in\N},\{\theta_k\}_{k\in\N}\subset\mathbb X$, and
	$\{y_k'\}_{k\in\N},\{\xi_k\}_{k\in\N}\subset\mathbb Y$
	which satisfy
	\begin{equation*}
		\begin{aligned}
		\max(\norm{x_k-\bar x_k},\nnorm{x_k'-\bar x_k})
		&<
		\min\left(\frac1k\norm{\bar x_k-\bar x},\dist(\bar x_k,\Phi^{-1}(\bar y))\right),
		\\
		\max(\norm{y_k-\bar y_k},\nnorm{y_k'-\bar y_k})
		&<
		\frac1k\norm{\bar y_k-\bar y},
		\\
		\norm{\theta_k}&<\norm{\bar x_k-\bar x},
		\\ 
		\norm{\xi_k}&<\frac1k\norm{\bar x_k-\bar x},
		\end{aligned}
	\end{equation*}		
	as well as
	\begin{equation}\label{eq:asymptotic_necessary_condition}
		\forall k\in\N\colon\quad
		\bar x-x_k'+\theta_k\in\widehat\partial\varphi(x_k')+\widehat D^*\Phi(x_k,y_k)(k(y_k'-\bar y)+\xi_k).
	\end{equation}
	Setting $\eta_k:=\bar x-x_k'+\theta_k$ for each $k\in\N$, we find $\eta_k\to 0$.
	The convergences $x_k,x_k'\to\bar x$, $y_k,y_k'\to\bar y$, $\theta_k\to 0$, and $\xi_k\to 0$ are easily obtained.
	For later use, we want to point out the estimates
	\begin{align*}
		\min(\norm{x_k-\bar x},\nnorm{x_k'-\bar x})&>\frac{k-1}k\norm{\bar x_k-\bar x},&
		\\
		\min(\norm{y_k-\bar y},\nnorm{y_k'-\bar y})&>\frac{k-1}k\norm{\bar y_k-\bar y}&	
	\end{align*}
	which follow from the reverse triangle inequality.
	They particularly yield $x_k,x_k'\neq \bar x$ and $y_k,y_k'\neq\bar y$ for each $k\in\N$.
	Furthermore, by $\bar x_k\notin\Phi^{-1}(\bar y)$, 
	$\norm{x_k-\bar x_k}<\dist(\bar x_k,\Phi^{-1}(\bar y))$, and 
	$\nnorm{x_k'-\bar x_k}<\dist(\bar x_k,\Phi^{-1}(\bar y))$,
	we already find $x_k,x_k'\notin\Phi^{-1}(\bar y)$ for each $k\in\N$.
	Since $\{(x_k - \bar x)/\nnorm{x_k - \bar x}\}_{k\in\N}\subset\mathbb S_{\mathbb X}$, 
	we may assume $(x_k - \bar x)/\nnorm{x_k - \bar x} \to u$ 
	for some $u\in\mathbb S_{\mathbb X}$.
	We introduce a sequence
	$\{y_k^*\}_{k\in\N}\subset\mathbb Y$ by means of 
	$y_k^* := (k(y_k'-\bar y)+\xi_k) \norm{y_k-\bar y}/\norm{x_k - \bar x}$ for each $k\in\N$.	
	
	First, we show that $(x_k'-\bar x)/\nnorm{x_k'-\bar x}, (\bar x_k-\bar x)/\nnorm{\bar x_k-\bar x}\to u$. This follows from
	\begin{equation}\label{dir_u}
		\norm{\frac{x_k'-\bar x}{\nnorm{x_k'-\bar x}}-\frac{x_k-\bar x}{\norm{x_k-\bar x}}}
		\leq
		\norm{\frac{x_k'-\bar x}{\nnorm{x_k'-\bar x}}-\frac{\bar x_k-\bar x}{\norm{\bar x_k-\bar x}}}
		+
		\norm{\frac{x_k-\bar x}{\norm{x_k-\bar x}}-\frac{\bar x_k-\bar x}{\norm{\bar x_k-\bar x}}}
	\end{equation}
	together with
	\begin{align*}
				&\norm{\frac{x_k-\bar x}{\norm{x_k-\bar x}}-\frac{\bar x_k-\bar x}{\norm{\bar x_k-\bar x}}}
				=
				\norm{\frac{x_k-\bar x_k}{\norm{x_k-\bar x}} + (\bar x_k-\bar x)\left(\frac{1}{\norm{x_k-\bar x}} - \frac{1}{\norm{\bar x_k-\bar x}}\right)}
				\\
				&\qquad
				\leq
				\frac{\norm{x_k-\bar x_k}}{\norm{x_k-\bar x}} + \frac{\norm{\bar x_k-\bar x}\norm{x_k-\bar x_k}}{\norm{x_k-\bar x}\norm{\bar x_k-\bar x}}	
				=	
				2\frac{\norm{x_k-\bar x_k}}{\norm{x_k-\bar x}}		
				\leq
				2\frac{\frac1k\norm{\bar x_k-\bar x}}{\frac{k-1}{k}\norm{\bar x_k-\bar x}}
				=
				\frac{2}{k-1},
	\end{align*}
	and the fact that analogous arguments show that the first summand on the right-hand side 
	of \eqref{dir_u} tends to zero, too.
	
	Next, we claim that $\{y_k^*\}_{k\in\N}$ is bounded.
	Due to
	\begin{equation}\label{eq:product_with_xi_not_critical}
		\frac{\nnorm{\xi_k}}{\nnorm{x_k - \bar x}}
		\leq
		\frac{1}{k-1}\frac{k\nnorm{\xi_k}}{\nnorm{\bar x_k - \bar x}}
		\leq \frac{1}{k-1}
	\end{equation}
	for all $k\in\N$, it suffices to show boundedness of $\{k (y_k'-\bar y)\nnorm{y_k-\bar y}/\nnorm{x_k-\bar x}\}_{k\in\N}$.
	First, we get the estimate
	\begin{equation}\label{al:bound}
		\begin{aligned}
		k \frac{\nnorm{y_k'-\bar y}\nnorm{y_k-\bar y}}{\nnorm{x_k-\bar x}}
		&\leq
		k\frac{(\nnorm{y_k'-\bar y_k}+\nnorm{\bar y_k-\bar y})(\nnorm{y_k-\bar y_k}+\nnorm{\bar y_k-\bar y})}{\frac{k-1}{k}\nnorm{\bar x_k-\bar x}}
		\\
		&
		\leq
		k\frac{(1 + \frac{1}{k})^2\nnorm{\bar y_k-\bar y}^2}{(1 - \frac{1}{k})\nnorm{\bar x_k-\bar x}}
		\leq
		2k\frac{\nnorm{\bar y_k-\bar y}^2}{\nnorm{\bar x_k-\bar x}}
		\end{aligned}
	\end{equation}
	for large enough $k\in\N$.
	Second, rearranging \eqref{eq:estimate_from_penalization}, leaving a nonnegative term away,
	and division by $\nnorm{\bar x_k-\bar x}$ give us
	\begin{equation}\label{eq:another_estimate_from_penalization}
		\forall k\in\N\colon
		\quad
		\frac{\varphi(\bar x_k)-\varphi(\bar x)}{\nnorm{\bar x_k-\bar x}}
		+
		\frac k2\frac{\nnorm{\bar y_k-\bar y}^2}{\nnorm{\bar x_k-\bar x}}
		\leq 
		0
		.
	\end{equation}
	Lipschitzianity of $\varphi$ yields boundedness of the first fraction, 
	so $\{k\nnorm{\bar y_k-\bar y}^2/\nnorm{\bar x_k-\bar x}\}_{k\in\N}$ must be bounded
	and, consequently, $\{y_k^*\}_{k\in\N}$ as well 
	by \eqref{eq:product_with_xi_not_critical} and \eqref{al:bound}.
	Thus, we may assume $y_k^*\to y^*$ for some $y^*\in\mathbb Y$.
	
	Suppose that $\{(y_k-\bar y)/\nnorm{x_k-\bar x}\}_{k\in\N}$ does not converge to
	zero. This, along a subsequence (without relabeling), yields boundedness
	of the sequence $\{k(y_k'-\bar y)+\xi_k\}_{k\in\N}$, 
	and taking the limit $k\to\infty$ in \eqref{eq:asymptotic_necessary_condition}
	yields~\ref{item:trivial_MSt}.
	
	Thus, we may assume $(y_k-\bar y)/\norm{x_k-\bar x}\to 0$ 
	for the reminder of the proof.
	Since we have $\nnorm{x_k-\bar x}\leq(1+1/k)\nnorm{\bar x_k-\bar x}$
	from the triangle inequality for each $k\in\N$, we also find
	\begin{equation}\label{eq:some_intermediate_convergence}
		\frac{\nnorm{\bar y_k-\bar y}}{\nnorm{\bar x_k-\bar x}}
		\leq
		\frac{\big(1+\frac{1}{k-1}\big)\nnorm{y_k-\bar y}}{\big(1-\frac{1}{k+1}\big)\nnorm{x_k-\bar x}}
		\to 0
	\end{equation}
	in this situation.
	Furthermore, observe that we have
	\begin{equation*}
		(\bar x+ \norm{\bar x_k - \bar x} (\bar x_k - \bar x)/\norm{\bar x_k - \bar x},\bar y+\norm{\bar x_k - \bar x} (\bar y_k - \bar y)/\norm{\bar x_k - \bar x})
		\in
		\gph\Phi
	\end{equation*}
	for all $k\in\N$.
	Additionally, \eqref{eq:another_estimate_from_penalization} yields 
	\begin{equation*}
		\limsup_{k\to\infty}
		\frac{\varphi(\bar x+ \norm{\bar x_k - \bar x} (\bar x_k - \bar x)/\norm{\bar x_k - \bar x})-\varphi(\bar x)}{\norm{\bar x_k - \bar x}}\leq 0,
	\end{equation*}
	so 
	\eqref{eq:some_intermediate_convergence} shows that 
	$u$ is a critical direction for \eqref{eq:nonsmooth_problem} at $\bar x$.
	
	In the remainder of the proof, we are going to exploit the sequences
	$\{v_k^\gamma\}_{k\in\N},\{\lambda_k^\gamma\}_{k\in\N}\subset\mathbb Y$ 
	given as in \eqref{eq:definition_of_surrogate_sequences}. 
	Observe that $y_k^* = \nnorm{v^\gamma_k} \lambda_k^\gamma $, 
	i.e., $\lambda_k^\gamma = y_k^*\norm{x_k - \bar x}^\gamma/\norm{y_k - \bar y}$ is valid for 
	each $k\in\N$.
	Note that the optimality conditions \eqref{eq:asymptotic_necessary_condition} can be rewritten as
	\begin{equation}\label{eq:asymptotic_necessary_condition_of_order_gamma}
		\forall k\in\N\colon\quad
		\eta_k
		\in
		\widehat\partial\varphi(x_k')+\widehat D^*\Phi(x_k,\yb 
		+
		 \norm{x_k - \bar x}^\gamma v^\gamma_k)\left(\frac{\lambda_k^\gamma}{\norm{x_k - \bar x}^{\gamma-1}}\right).
	\end{equation}	
	Now, we need to distinguish three options.

	Let us assume that $\lambda_k^\gamma\to 0$.
	Using $t_k:=\norm{x_k-\bar x}$, 
	we can reformulate \eqref{eq:asymptotic_necessary_condition_of_order_gamma} as
	\[
		\forall k\in\N\colon\quad
		\eta_k
		\in
		\widehat\partial\varphi(x_k')
		+
		\widehat D^*\Phi\left(\bar x+t_k\frac{x_k-\bar x}{\norm{x_k-\bar x}},\bar y+t_k\frac{y_k-\bar y}{\norm{x_k-\bar x}}\right)
		\left(\frac{\lambda_k^\gamma}{t_k^{\gamma-1}}\right).
	\]
	Taking the limit $k\to\infty$ yields~\ref{item:pseudo_MSt_zero_mult}
	since $(y_k-\bar y)/\norm{x_k-\bar x}\to 0$ 
	and $u$ has already been shown to be critical for \eqref{eq:nonsmooth_problem} at $\bar x$.

	If $\{\lambda_k^\gamma\}_{k\in\N}$ remains bounded but, along a subsequence
	(without relabeling), stays away from zero, 
	we also get boundedness of $\{v_k^{\gamma}\}_{k\in\N}$
	from boundedness of $\{y_k^*\}_{k\in\N}$,
	and taking the limit along a convergent subsequence (without relabeling) in 
	\eqref{eq:asymptotic_necessary_condition_of_order_gamma} yields precisely \eqref{eq:Pseudo-M-stat}, 
	where $v,\lambda\in\mathbb Y$ 
	with $\lambda\neq 0$
	satisfy
	$v_k^{\gamma}\to v$ and $\lambda_k^{\gamma}\to\lambda$, respectively.
	The additional information $\innerprod{\lambda}{ v }\geq 0$ comes from the estimate
	\begin{align*}
		\innerprod{y_k' - \bar y}{y_k - \bar y}
		& = 
		\norm{y_k - \bar y}^2 + \innerprod{y_k' - \bar y_k + \bar y_k - y_k}{y_k - \bar y}\\
		& \geq 
		\norm{y_k - \bar y}^2 - \norm{y_k - \bar y}\big(\nnorm{y_k' - \bar y_k} + \norm{\bar y_k - y_k} \big)\\
		& \geq  
		\norm{y_k - \bar y}^2\big(1 - 2/k\big),
	\end{align*}
	which yields
	\begin{align*}
		\innerprod{\lambda}{ v }
		& = 
		\lim_{k \to \infty} \innerprod{\lambda_k^\gamma }{v^\gamma_k }
		 = 
		\lim_{k \to \infty} \left( \frac{k}{\norm{x_k - \bar x}} 
		\innerprod{y_k' - \bar y}{y_k - \bar y} 
		+  
		\frac{1}{\norm{x_k - \bar x}}\innerprod{\xi_k}{y_k - \bar y} \right)\\
		& \geq 
		\lim_{k \to \infty} - \nnorm{\xi_k}  \frac{\norm{y_k - \bar y}}{\norm{x_k - \bar x}}
		= 
		0
	\end{align*}
	where the final limit follows from \eqref{eq:product_with_xi_not_critical}.
	Criticality of $u$ for \eqref{eq:nonsmooth_problem} at $\bar x$ 
	has been shown above.

	Finally, if $\{\lambda_k^\gamma\}_{k\in\N}$ is not bounded, we pass to a subsequence (without relabeling) such that $\nnorm{\lambda_k^{\gamma}} \to \infty$
	and so we also get $v_k^{\gamma} \to 0$ along this subsequence by boundedness of $\{y_k^*\}_{k\in\N}$.
	It remains to show $(y_k-\bar y)/\norm{y_k-\bar y}-\lambda_k^\gamma/\nnorm{\lambda_k^\gamma} \to 0$
	in order to validate~\ref{item:irregular_case_pseudo_MSt}.
	We have
	\begin{align*}
		&\norm{\frac{y_k-\bar y}{\norm{y_k-\bar y}}-\frac{y_k'-\bar y}{\nnorm{y_k'-\bar y}}}
		=
		\frac{\norm{\nnorm{y_k'-\bar y}(y_k-\bar y)-\norm{y_k-\bar y}(y_k'-\bar y)}}{\norm{y_k-\bar y}{\nnorm{y_k'-\bar y}}}\\
		&\qquad
		=
		\frac{\norm{\nnorm{y_k'-\bar y}[(y_k-\bar y)-(y_k'-\bar y)]+[\nnorm{y_k'-\bar y}-\norm{y_k-\bar y}](y_k'-\bar y)}}{\norm{y_k-\bar y}{\nnorm{y_k'-\bar y}}}\\
		&\qquad
		\leq
		2\frac{\nnorm{y_k-y_k'}}{\norm{y_k-\bar y}}
		\leq
		2\frac{\norm{y_k-\bar y_k}+\nnorm{\bar y_k-y_k'}}{\norm{y_k-\bar y}}
		\leq
		2\frac{\frac2k\norm{\bar y_k-\bar y}}{\frac{k-1}k\norm{\bar y_k-\bar y}}
		=
		\frac{4}{k-1}\to 0.
	\end{align*}
	On the other hand, since $\xi_k \to 0$ and 
	\[
		k\nnorm{y_k'-\bar y} 
		\geq
		\norm{y_k^*}\frac{\norm{x_k-\bar x}}{\norm{y_k-\bar y}}
		- 
		\norm{\xi_k} 
		=
		\nnorm{\lambda_k^\gamma}\norm{x_k-\bar x}^{1-\gamma} 
		-
		\norm{\xi_k}
		\to 
		\infty,
	\]
	we get $k\norm{y_k'-\bar y} \to \infty$ and thus
	\begin{align*}
		\lim\limits_{k\to\infty}\frac{y_k-\bar y}{\norm{y_k-\bar y}}
		&=
		\lim\limits_{k\to\infty}\frac{y_k'-\bar y}{\nnorm{y_k'-\bar y}}
		=
		\lim\limits_{k\to\infty}\frac{k(y_k'-\bar y)+\xi_k}{k\nnorm{y_k'-\bar y}}
		\\
		&
		=
		\lim\limits_{k\to\infty}
			\frac{k(y_k'-\bar y)+\xi_k}{\nnorm{k\nnorm{y_k'-\bar y}+\xi_k}}
			\frac{\nnorm{k\nnorm{y_k'-\bar y}+\xi_k}}{k\nnorm{y_k'-\bar y}}
		\\
		&
		=
		\lim\limits_{k\to\infty}\frac{y_k^*}{\nnorm{y_k^*}}
		=
		\lim\limits_{k\to\infty}\frac{\lambda_k^\gamma}{\nnorm{\lambda_k^\gamma}}.
	\end{align*}
	Keeping in mind that $u$ is critical for \eqref{eq:nonsmooth_problem} at $\bar x$,
	all conditions stated in~\ref{item:irregular_case_pseudo_MSt} have been verified
	since \eqref{eq:asymptotic_stationarity_regular_tools_gamma} follows from
	\eqref{eq:asymptotic_necessary_condition_of_order_gamma}.
	
	Finally, let us argue that option~\ref{item:irregular_case_pseudo_MSt} can be avoided, i.e., that the sequence
	$\{\lambda_k^\gamma\}_{k\in\N}$ from above remains bounded,
	if we assume that $\Phi$ is metrically pseudo-subregular
	of order $\gamma$ in direction $u$ at $(\xb,\yb)$.
	From \eqref{eq:some_intermediate_convergence}, we find
	\begin{equation*}
		\frac{\norm{\bar y_k-\bar y}}{\norm{\bar x_k-\bar x}^\gamma}
		\leq
		\frac{\big(1 + \frac{1}{k-1}\big)}{\big(1 - \frac{1}{k+1}\big)^\gamma} \frac{\norm{y_k-\bar y}}{\norm{x_k-\bar x}^\gamma} 
		\leq
		2^\gamma
		\nnorm{v_k^{\gamma}} + \oo(\nnorm{v_k^{\gamma}})
	\end{equation*}
	for each $k\in\N$.
	Thus, by metric pseudo-subregularity of $\Phi$, 
	there is a constant $\kappa>0$ such that, for sufficiently large $k\in\N$, we get the existence of $\tilde x_k \in \Phi^{-1}(\yb)$ with
	\begin{equation}\label{eq:estimate_pseudo_subregularity}
		\norm{\bar x_k - \tilde x_k} 
		\leq
		\kappa \frac{\dist(\bar y,\Phi(\bar x_k))}{\norm{\bar x_k-\bar x}^{\gamma - 1}} 
		\leq
		\kappa \frac{\norm{\bar y_k - \bar y}}{\norm{\bar x_k-\bar x}^{\gamma - 1}} 
		\leq
		\kappa \norm{\bar x_k-\bar x} \bigl(2^\gamma\nnorm{v_k^{\gamma}} + \oo(\nnorm{v_k^{\gamma}})\bigr).
	\end{equation}
	Since $(\bar x_k,\bar y_k)$ 
	is a global minimizer of \eqref{eq:penalized_nonsmooth_problem}, we get
	\[
		\frac{\varphi(\bar x_k)-\varphi(\tilde x_k)}{\norm{\bar x_k-\bar x}}
		+
		\frac k2\frac{\norm{\bar y_k-\bar y}^2}{\norm{\bar x_k-\bar x}}
		+
		\frac 12\frac{\norm{\bar x_k-\bar x}^2 - \norm{\tilde x_k-\bar x}^2}{\norm{\bar x_k-\bar x}}
		\leq 
		0
	\]
	for all sufficiently large $k\in\N$.
	By boundedness of $\{y_k^*\}_{k\in\N}$, we immediately obtain the boundedness of $\{\lambda_k^\gamma\}_{k\in\N}$
	unless we have $v_k^\gamma\to 0$. Thus, let us assume the latter.
	Particularly, we find $\norm{\bar x_k - \tilde x_k} \to 0$ from above.
	Thus, taking into account $\xi_k \to 0$, \eqref{al:bound}, the above estimate, Lipschitzianity of $\varphi$, 
	and \eqref{eq:estimate_pseudo_subregularity}, we conclude
	\begin{align*}
        \nnorm{\lambda_k^{\gamma}}\nnorm{v_k^{\gamma}}
		&\leq
		2k\frac{\norm{\bar y_k-\bar y}^2}{\norm{\bar x_k-\bar x}} + \oo(\nnorm{v_k^{\gamma}})
		\\
		&\leq
		4L \frac{\norm{\bar x_k-\tilde x_k}}{\norm{\bar x_k-\bar x}} 
		+ 
		2\frac{(\norm{\tilde x_k-\bar x} - \norm{\bar x_k-\bar x})(\norm{\tilde x_k-\bar x} + \norm{\bar x_k-\bar x})}{\norm{\bar x_k-\bar x}} 
		+ 
		\oo(\nnorm{v_k^{\gamma}})
		\\
		&\leq
		\frac{\norm{\bar x_k-\tilde x_k}}{\norm{\bar x_k-\bar x}} \left( 4L + 2 (\norm{\bar x_k-\bar x} + \norm{\tilde x_k-\bar x})\right) 
		+ 
		\oo(\nnorm{v_k^{\gamma}})
		\\
		&\leq
		\kappa (2^\gamma\nnorm{v_k^{\gamma}} + \oo(\nnorm{v_k^{\gamma}})) \left( 4L + \oo(1)\right) + \oo(\nnorm{v_k^{\gamma}})
		\\
		&
		= 
		2^{2+\gamma}L \kappa  \nnorm{v_k^{\gamma}} + \oo(\nnorm{v_k^{\gamma}})
	\end{align*}
	for some constant $L>0$ and large enough $k\in\N$.
	Consequently, we find
	\[
		\nnorm{\lambda_k^{\gamma}} \leq 2^{2+\gamma}L \kappa  + 1
	\]
	for large enough $k\in\N$,
	and boundedness of $\{\lambda_k^\gamma\}_{k\in\N}$ follows.
\end{proof}

In case where $\gamma:=1$ holds true, both notions of a pseudo-coderivative from \eqref{def:coderivatives}
coincide with the directional limiting coderivative which is why we obtain the following
corollary, which enhances \cite[Theorem~3.9]{Mehlitz2020a}, from \cref{thm:higher_order_directional_asymptotic_stationarity}.
\begin{corollary}\label{thm:directional_asymptotic_stationarity}
	Let $\bar x\in\mathcal F$ be a local minimizer of \eqref{eq:nonsmooth_problem}.
	Then $\bar x$ is M-stationary or there exist a critical direction 
	$u\in\mathbb S_{\mathbb X}$ for \eqref{eq:nonsmooth_problem} at $\bar x$
	and $y^* \in \mathbb Y$ 
	as well as sequences
	$\{x_k\}_{k\in\N},\{x_k'\}_{k\in\N},\{\eta_k\}_{k\in\N}\subset\mathbb X$ and
	$\{y_k\}_{k\in\N},\{y_k^*\}_{k\in\N}\subset\mathbb Y$ such that
	$x_k,x_k'\notin\Phi^{-1}(\bar y)$, $y_k\neq\bar y$, and $y_k^*\neq 0$ for all $k\in\N$,
	\begin{equation}\label{eq:convergences_gamma=1}
		\begin{aligned}
			x_k&\to\bar x,&		\qquad	x_k'&\to\bar x,&	\qquad	y_k&\to\bar y,&
			\\
			\frac{x_k-\bar x}{\nnorm{x_k-\bar x}}&\to u,&	\qquad
			\frac{x_k'-\bar x}{\nnorm{x_k'-\bar x}}&\to u,&	\qquad
			\frac{y_k-\bar y}{\nnorm{x_k-\bar x}}&\to 0,&
			\\
			y_k^*&\to y^*,&\qquad
			\frac{\norm{x_k - \bar x}}{\norm{y_k - \bar y}}\norm{y_k^*} &\to \infty,&\qquad
			\eta_k&\to 0,&
			\\
			\frac{y_k-\bar y}{\nnorm{y_k-\bar y}}-\frac{y_k^*}{\nnorm{y_k^*}}&\to 0,& &&&&
		\end{aligned}
	\end{equation}
	and
	\begin{equation}\label{eq:asymptotic_stationarity_gamma=1}
		\forall k\in\N\colon\quad
		\eta_k\in\widehat{\partial}\varphi(x_k')+\widehat{D}^*\Phi(x_k,y_k)\left(\frac{\norm{x_k - \bar x}}{\norm{y_k - \bar y}}y_k^*\right).
	\end{equation}
\end{corollary}

The above result shows that each local minimizer of \eqref{eq:nonsmooth_problem}
either is M-stationary or satisfies so-called approximate (or asymptotic)
stationarity conditions w.r.t.\ a certain critical direction and an unbounded
sequence of multipliers. Related results in non-directional form can be found
in \cite{KrugerMehlitz2021,Mehlitz2020a}. The story of approximate stationarity
conditions in variational analysis, however, can be traced back to
\cite{Kruger1985,KrugerMordukhovich1980}. This concept has been rediscovered as a valuable
tool for the analysis of convergence properties for solution algorithms associated
with standard nonlinear optimization problems in 
\cite{AndreaniMartinezSvaiter2010,AndreaniHaeserMartinez2011}, and extensions were
made to disjunctive, conic, and even infinite-dimensional optimization,
see e.g.\ \cite{AndreaniHaeserSecchinSilva2019,AndreaniGomezHaeserMitoRamos2021,BoergensKanzowMehlitzWachsmuth2019,Ramos2019}
and the references therein.
Taking a look back at \cref{thm:higher_order_directional_asymptotic_stationarity},
one might be tempted to refer to the situation described 
in~\ref{item:irregular_case_pseudo_MSt}
as directional approximate stationarity of order $\gamma$.
\begin{remark}\label{rem:approximate_stationarity_wrt_all_critical_directions}
	It is well known from \cite[Theorem~7(ii)]{Gfrerer2013} 
	that under suitable constraint qualifications,
	M-stationarity-type stationarity conditions hold at a given local 
	minimizer w.r.t.\ \emph{each}
	critical direction. 
	Our result from \cref{thm:directional_asymptotic_stationarity} is essentially
	different since it shows that each local minimizer is M-stationary or satisfies approximate
	M-stationarity-type conditions for \emph{one} critical direction. It remains open whether
	the assertion of \cref{thm:directional_asymptotic_stationarity} can be extended to hold for
	\emph{each} critical direction since, unfortunately, we also do not have a counterexample 
	available which shows that this is not true.
\end{remark}

As a consequence of \cref{thm:directional_asymptotic_stationarity}, we obtain the following 
FJ-type necessary optimality condition for \eqref{eq:nonsmooth_problem}, 
see \cite[Lemma~3.4]{Mehlitz2020a} for a related result.
\begin{corollary}\label{cor:Fritz_John_type}
	Let $\bar x\in\mathcal F$ be a local minimizer of \eqref{eq:nonsmooth_problem}.
	Then $\bar x$ is M-stationary or there exist a critical direction 
	$u\in\mathbb S_{\mathbb X}$
	for \eqref{eq:nonsmooth_problem} at $\bar x$ and a multiplier 
	$\lambda\in\mathbb Y\setminus\{0\}$
	such that $0\in D^*\Phi((\bar x,\bar y);(u,0))(\lambda)$.
\end{corollary}
\begin{proof}
	Supposing that $\bar x$ is not M-stationary, \cref{thm:directional_asymptotic_stationarity}
	yields the existence of sequences $\{x_k\}_{k\in\N},\{x_k'\}_{k\in\N},\{\eta_k\}_{k\in\N}\subset\mathbb X$
	and $\{y_k\}_{k\in\N},\{y_k^*\}_{k\in\N}\subset\mathbb Y$ satisfying $x_k,x_k'\neq\bar x$,
	$y_k\neq\bar y$, and $y_k^*\neq 0$ for all $k\in\N$ as well as \eqref{eq:convergences_gamma=1}
	and \eqref{eq:asymptotic_stationarity_gamma=1}.
	Setting $\lambda_k:=y_k^*\norm{x_k-\bar x}/\norm{y_k-\bar y}$ for each $k\in\N$, we find $\norm{\lambda_k}\to\infty$
	from \eqref{eq:convergences_gamma=1}. 
	Using this notation, for each $k\in\N$, we find $x_k^*\in\widehat{\partial}\varphi(x_k')$ such that $\eta_k-x_k^*\in\widehat D^*\Phi(x_k,y_k)(\lambda_k)$.
	Dividing the latter by $\norm{\lambda_k}$ while exploiting positive homogeneity of the
	regular coderivative and boundedness of $\{x_k^*\}_{k\in\N}$ which follows by Lipschitzianity of $\varphi$,
	taking the limit along a suitable subsequence while respecting the definition of the directional limiting
	coderivative gives us some nonvanishing $\lambda\in\mathbb Y$ such that $0\in D^*\Phi((\bar x,\bar y);(u,0))(\lambda)$
	holds.
\end{proof}

Clearly, \cref{cor:Fritz_John_type} illustrates nicely that FOSCMS
provides a constraint qualification implying the local minimizer $\bar x$ of \eqref{eq:nonsmooth_problem}
to be M-stationary. As mentioned earlier, this already follows from \cite[Theorem~7(ii)]{Gfrerer2013}.
However, \cref{thm:directional_asymptotic_stationarity} can also be used to formulate so-called
\emph{asymptotic} constraint qualifications to rule out the situation where M-stationarity fails.
In the non-directional situation, this has been done in \cite{Mehlitz2020a}.
Part B of this paper is devoted to an associated theory in the directional situation.

Let us now come back to our necessary optimality conditions from \cref{thm:higher_order_directional_asymptotic_stationarity}.
In the subsequently stated remark, we comment on some additional information \cref{thm:higher_order_directional_asymptotic_stationarity}
provides about the involved subgradients of the objective function $\varphi$.

\begin{remark}\label{rem:directional_subgradients_in_stationarity_condition}
	The proof of \cref{thm:higher_order_directional_asymptotic_stationarity}
	shows that the sequence $\{x_k'\}_{k\in\N}$ additionally satisfies
	\begin{align*}
		- L \nnorm{x_k' - \bar x}
		\leq 
		\varphi(x_k') - \varphi(\bar x) 
		& \leq 
		\varphi(\bar x_k) - \varphi(\bar x) + \varphi(x_k') - \varphi(\bar x_k) 
		\leq 
		L \nnorm{x_k' - \bar x_k}
		\\
		& \leq 
		L/k \norm{\bar x_k - \bar x}
		<  
		L/(k-1) \nnorm{x_k' - \bar x},
	\end{align*}
	taking into account $\varphi(\bar x_k) - \varphi(\bar x) \leq 0$ for each $k\in\N$, which follows 
	from \eqref{eq:estimate_from_penalization}.
	Here, $L>0$ is some local Lipschitz constant of $\varphi$ valid around $\bar x$.
	Thus, by passing to a subsequence (without relabeling) we may assume that
	\[\frac{\varphi(x_k') - \varphi(\bar x)}{\nnorm{x_k' - \bar x}} \to \mu \leq 0.\]
	This can be added to case~\ref{item:irregular_case_pseudo_MSt} of \cref{thm:higher_order_directional_asymptotic_stationarity}.
	
	Moreover, since, for each $k\in\N$, $x_k^* \in \widehat\partial \varphi(x_k')$ means
	\[
		(x_k^*,-1)
		\in
		\widehat{\mathcal N}_{\epi \varphi}(x_k',\varphi(x_k'))
		=
		\widehat{\mathcal N}_{\epi \varphi}(\xb + t_k (x_k' - \xb)/t_k,\varphi(\xb) + t_k (\varphi(x_k') - \varphi(\xb))/t_k)
	\]
	for $t_k := \nnorm{x_k' - \bar x}$, for the limit $x^*\in\mathbb X$ of $\{x_k^*\}_{k\in\N}$, we get
	\[
		x^*
		\in 
		\partial \varphi(\bar x;(u,\mu)) 
		:= 
		\{x^*\in\mathbb X \,|\, 
			(x^*,-1) \in \mathcal N_{\epi \varphi}((\bar x,\varphi(\bar x));(u,\mu))
		\},
	\]
	where $\partial \varphi(\bar x;(u,\mu))$
	denotes the (geometric) subdifferential of $\varphi$ at $\bar x$ in direction $(u,\mu)$, see \cite{BenkoGfrererOutrata2019}.
	Consequently, in cases~\ref{item:pseudo_MSt_zero_mult} and~\ref{item:pseudo_dir_MSt}
	of \cref{thm:higher_order_directional_asymptotic_stationarity}, $\partial \varphi(\bar x)$ can be replaced by $\partial \varphi(\bar x;(u,\mu))$.
\end{remark}

From \cref{lem:sufficient_condition_pseudo_subregularity_abstract} and \cref{thm:higher_order_directional_asymptotic_stationarity}, 
we get the following corollary.

\begin{corollary}\label{Cor:M_stat_order_gamma}
	Let $\bar x\in\mathcal F$ be a local minimizer of \eqref{eq:nonsmooth_problem} and consider $\gamma \geq 1$.
	Assume that \eqref{eq:FOSCMS_gamma} holds for each 
	$u\in\ker D\Phi(\bar x,\bar y)\cap \mathbb S_{\mathbb X}$.
    Then one of the alternatives~\ref{item:trivial_MSt},~\ref{item:pseudo_MSt_zero_mult}, or~\ref{item:pseudo_dir_MSt} 
    from \cref{thm:higher_order_directional_asymptotic_stationarity} is valid.
\end{corollary}

Observe that \cref{thm:higher_order_directional_asymptotic_stationarity} and \cref{Cor:M_stat_order_gamma} 
foreshadow a new approach on how to identify conditions which
guarantee that local minimizers $\bar x\in\mathcal F$ of \eqref{eq:nonsmooth_problem} are M-stationary.
First, we need to fix some $\gamma>1$.
Next, we employ \cref{lem:sufficient_condition_pseudo_subregularity_abstract}, particularly, condition \eqref{eq:FOSCMS_gamma}
which is stated  in terms of the pseudo-coderivative, in order to rule out option~\ref{item:irregular_case_pseudo_MSt} 
of \cref{thm:higher_order_directional_asymptotic_stationarity}.
Now, we know that $\bar x$ is either M-stationary, or
one of the stationarity conditions \eqref{eq:Pseudo-M-Stat_Gferer} and \eqref{eq:Pseudo-M-stat} 
from cases~\ref{item:pseudo_MSt_zero_mult} and~\ref{item:pseudo_dir_MSt},
which can be written in a unified way as
\begin{equation}\label{eq:unified}	
	0 \in \partial\varphi(\bar x) + \widetilde{D}^\ast_{\gamma} \Phi((\xb,\yb);(u,0))(\lambda),
\end{equation}
see \eqref{eq:trivial_upper_estimate_pseudo_coderivative} as well, is valid.
Finally, we need to check whether the image of Gfrerer's directional pseudo-coderivative is 
included in the image of the standard limiting coderivative.
This procedure directly leads to the problem of how to compute or estimate the appearing pseudo-coderivatives.
Exemplary, for $\gamma:=2$ and in the setting where $\Phi$ is a constraint mapping, these objects actually can be
computed, see \cref{sec:variational_analysis_constraint_mapping}, and we obtain conditions in terms of initial problem data. 
This is discussed in detail in \cref{sec:constraint_mappings}.

\subsection{New necessary optimality conditions and constraint qualifications for optimization problems with geometric constraints}
\label{sec:constraint_mappings}

Again, we are concerned with necessary optimality conditions for problem
\eqref{eq:nonsmooth_problem}.
Here, we look more closely into the case where 
$\Phi\colon\mathbb X\tto\mathbb Y$ is given in
the form of a constraint mapping, i.e., $\Phi(x) := g(x) - D$ for all $x\in\mathbb X$ 
may hold where $g\colon\mathbb X\to\mathbb Y$ is twice continuously differentiable 
and $D \subset \mathbb Y$ is a closed set.
Furthermore, we focus on the order $\gamma:=2$ in our considerations, so that we can exploit
our results from \cref{sec:variational_analysis_constraint_mapping} 
where we computed the pseudo-coderivative of order $2$ of $\Phi$.
For that purpose, we assume $\bar y:=0$ in \eqref{eq:nonsmooth_problem} 
throughout the section which can be done
without loss of generality as already pointed out in \cref{sec:directional_asymptotic_tools}.
Based on \cref{The : NCgen}, we obtain the following result.

\begin{proposition}\label{Pro:M-stat_via_second_order}
 Let $\bar x\in\mathcal F$ be a local minimizer of \eqref{eq:nonsmooth_problem}.
 	\begin{enumerate}
  	\item\label{item:pseudo_stationarity_constraint_maps} 
  		If, for every $u \in \mathbb S_{\mathbb X}$, we have \eqref{eq:CQ_pseudo_subregularity_II},
  		as well as \eqref{eq:CQ_pseudo_subregularity_Ia} or, in case $\nabla g(\bar x)u\neq 0$,
  		\eqref{eq:CQ_pseudo_subregularity_Ib},
   		then $\xb$ is M-stationary or there exist some critical direction 
   		$u\in\mathbb S_{\mathbb X}$ of \eqref{eq:nonsmooth_problem} at $\bar x$ and
  		\[
  			y^* \in \mathcal N_{D}(g(\xb);\nabla g(\xb)u)  \cap \ker\nabla g(\xb)^*,
  			\qquad 	
  			z^* \in D\mathcal N_{D}(g(\xb),y^*)(\nabla g(\xb)u)
  		\]
  		such that
  		\begin{equation}\label{eq:Sec_order_optim_cond}
   				0 
   				\in 
   				\partial\varphi(\bar x) 
   				+ 
   				\nabla^2\innerprod{y^*}{g}(\xb)(u) + \nabla g(\xb)^* z^*.
  		\end{equation}
  		If, moreover, for every such $u$, $y^*$, and $z^*$, 
  		there exists $\lambda \in \mathcal N_{D}(g(\xb))$ with
  		\begin{equation}\label{eq:Pseudo_to_M_stat}
    		\nabla^2\innerprod{y^*}{g}(\xb)(u) + \nabla g(\xb)^* z^* 
    		= 
    		\nabla g(\xb)^* \lambda,
  		\end{equation}
  		then $\xb$ is M-stationary.
  	\item\label{item:pseudo_stationarity_constraint_maps_polyhedral} 
  		If $\mathbb Y:=\R^m$, $D$ is locally polyhedral around $g(\xb)$, 
  		and condition \eqref{eq:CQ_pseudo_subregularity_polyhedral_I}
  		holds for every $u \in \mathbb S_{\mathbb X}$, 
  		then $\xb$ is M-stationary or there exist a critical direction 
  		$u\in\mathbb S_{\mathbb X}$ of \eqref{eq:nonsmooth_problem} at $\bar x$ and
  		\[
  			y^* \in \mathcal N_{\mathcal T_D(g(\bar x))}(\nabla g(\bar x) u) 
  			\cap 
  			\ker\nabla g(\xb)^*,
  			\qquad
  			z^* \in \mathcal T_{ \mathcal N_{\mathcal T_D(g(\bar x))}(\nabla g(\bar x) u)}(y^*)
  		\]
  		satisfying \eqref{eq:Sec_order_optim_cond}.
  		If, moreover, for every such $u$, $y^*$, and $z^*$, 
  		there exists $\lambda \in \mathcal N_{D}(g(\xb))$ with \eqref{eq:Pseudo_to_M_stat},
  		then $\xb$ is M-stationary.
  	\item\label{item:pseudo_stationarity_constraint_maps_polyhedral_refined} 
  		If $\mathbb Y:=\R^m$, $D$ is locally polyhedral around $g(\xb)$, 
  		and condition \eqref{eq:CQ_pseudo_subregularity_polyhedral_II}
  		holds for every $u \in \mathbb S_{\mathbb X}$ and $s\in\mathbb X$,
  		then $\xb$ is M-stationary or there exist a critical direction 
  		$u\in\mathbb S_{\mathbb X}$ of \eqref{eq:nonsmooth_problem} at $\bar x$, 
  		$s\in\mathbb X$, $v \in \mathbb Y$, and
  		\[
  			y^* \in \mathcal N_{\mathbf T(u)}(w_s(u,v)) \cap \ker\nabla g(\xb)^*,
  	 		\qquad
  	 		z^* \in \mathcal T_{\mathcal N_{\mathbf T(u)}(w_s(u,v))}(y^*)
  		\]
  		satisfying \eqref{eq:Sec_order_optim_cond} 
  		as well as $\innerprod{v}{ y^* }\geq 0$,
  		where $w_s(u,v)$ has been defined in \eqref{eq:Tu_and_ws}.
  		If, moreover, for every such $u$, $v$, $s$, $y^*$, and $z^*$, 
  		there exists $\lambda \in \mathcal N_{D}(g(\xb))$ with \eqref{eq:Pseudo_to_M_stat},
  		then $\xb$ is M-stationary.
 	\end{enumerate}
\end{proposition}
\begin{proof}
		For the proof of~\ref{item:pseudo_stationarity_constraint_maps},
		we first apply \cref{cor:sufficient_condition_pseudo_subregularity}
		in order to find that
		$\Phi$ is metrically pseudo-subregular of order $2$ at $(\bar x,\bar y)$ 
		in every direction $u\in\mathbb S_{\mathbb X}$.
		Then \cref{Cor:M_stat_order_gamma} and the subsequently stated discussions 
		imply that $\bar x$ is M-stationary or
		\eqref{eq:unified} holds for some $\lambda \in \mathbb Y$. 
		Utilizing \eqref{eq:CQ_pseudo_subregularity_Ia} or 
		$\nabla g(\bar x)u\neq 0$ and \eqref{eq:CQ_pseudo_subregularity_Ib} 
		together with \cref{The : NCgen}\,\ref{item:general_estimate_+CQ} 
		gives us $z^*\in\mathbb Y$
		with the required properties.
		The final claim about M-stationarity is obvious.

		The arguments in the polyhedral 
		case~\ref{item:pseudo_stationarity_constraint_maps_polyhedral} are analogous, 
		relying on \cref{The : NCgen}\,\ref{item:polyhedral_estimate} instead.

		For the proof of~\ref{item:pseudo_stationarity_constraint_maps_polyhedral_refined},
		similarly as in the previous cases, we find that $\Phi$ is metrically
		pseudo-subregular of order $2$ at $(\bar x,\bar y)$ 
		in every direction $u \in \mathbb S_{\mathbb X}$ due to 
		\cref{cor:sufficient_condition_pseudo_subregularity}.
		\Cref{Cor:M_stat_order_gamma} then implies that $\bar x$ is M-stationary or
		one of the cases~\ref{item:pseudo_MSt_zero_mult} and~\ref{item:pseudo_dir_MSt}
		from \cref{thm:higher_order_directional_asymptotic_stationarity} holds.
		In case of 
		\cref{thm:higher_order_directional_asymptotic_stationarity}\,\ref{item:pseudo_MSt_zero_mult}, 
		however, from \cref{The : NCgen}\,\ref{item:polyhedral_estimate} we get
		$0 \in \partial \varphi(\xb) + \nabla g(\xb)^* z^*$ 
		with $z^* \in \mathcal T_{ \mathcal N_{\mathcal T_D(g(\bar x))}(\nabla g(\bar x) u)}(0)$.
		This means 
		$z^* \in \mathcal N_{\mathcal T_D(g(\bar x))}(\nabla g(\bar x) u) 
		=\mathcal N_D(g(\bar x);\nabla g(\bar x)u)\subset \mathcal N_D(g(\bar x))$,
		see \cref{lem:some_properties_of_polyhedral_sets}, 
		and M-stationarity of $\bar x$ follows.
		In case of \cref{thm:higher_order_directional_asymptotic_stationarity}\,\ref{item:pseudo_dir_MSt}, 
		from \cref{The : NCgen}\,\ref{item:polyhedral_estimate_order_2} 
		we precisely obtain $y^*$ and $z^*$ as stated.
\end{proof}

Let us recall that although more technical, the assumptions in 
\cref{Pro:M-stat_via_second_order}\,\ref{item:pseudo_stationarity_constraint_maps_polyhedral_refined}
are milder than those ones in statement~\ref{item:pseudo_stationarity_constraint_maps_polyhedral},
see \cref{rem:upper_estimate_pseudo_coderivative_constraint_maps}.

The above proposition offers conditions that yield
M-stationarity in the general and in the polyhedral case.
In Part B of this paper, we recover the same results from a different angle,
namely, we show that these conditions are sufficient for so-called {\em directional asymptotic regularity}
(which then yields M-stationarity).
Moreover, we prove that these conditions are not stronger than the aforementioned
condition FOSCMS, see \cref{sec:generalized_differentiation}, and the
Second-Order Sufficient Condition for Metric Subregularity, see e.g.\ \cite{GfrererKlatte2016},
which also serve as constraint qualifications guaranteeing
M-stationarity of local minimizers in the present setting.

\section{Applications}\label{sec:applications}

In this section, we specify the results of \cref{sec:constraint_mappings} to some popular settings
in optimization theory. More precisely, we focus on the feasible regions of standard nonlinear,
complementarity-constrained, and nonlinear semidefinite problems.

\subsection{Standard nonlinear programming}\label{sec:NLPs}

We address a twice continuously differentiable function $g\colon\mathbb X\to\R^m$ 
with component functions
$g_1,\ldots,g_m\colon\mathbb X\to\R$, and, for some $\ell\in\{0,\ldots,m\}$, 
we study the constraint system
\begin{equation}\label{eq:NLP}\tag{NLP}
	g_i(x)\leq 0\quad i\in\{1,\ldots,\ell\},\qquad g_i(x)=0\quad i\in\{\ell+1,\ldots,m\}.
\end{equation}
Note that this corresponds to the investigation of the constraint map 
$\Phi\colon\mathbb X\tto\R^m$ given by
$\Phi(x):=g(x)-D$, $x\in\mathbb X$, 
where $D$ is the Cartesian product of $\R^{\ell}_-$ and the zero in $\R^{m-\ell}$, i.e.,
a convex, polyhedral cone.
For a given point $\bar x\in\mathbb X$ feasible to \eqref{eq:NLP}, we introduce
\[
	I(\bar x):=\{i\in\{1,\ldots,\ell\}\,|\,g_i(\bar x)=0\},\qquad
	I^-(\bar x):=\{i\in\{1,\ldots,\ell\}\,|\,g_i(\bar x)<0\}
\]
and note that each critical direction $u\in\mathbb S_{\mathbb X}$ of the associated problem \eqref{eq:nonsmooth_problem} must satisfy
\begin{equation}\label{eq:linearization_cone_NLP}
	\nabla g_i(\bar x) u\leq 0\quad i\in I(\bar x) ,\qquad 
	\nabla g_i(\bar x) u=0\quad i\in\{\ell+1,\ldots,m\}.
\end{equation}
Now, we can specify \eqref{eq:CQ_pseudo_subregularity_polyhedral_I}
which, for some $u\in\mathbb S_{\mathbb X}$ satisfying \eqref{eq:linearization_cone_NLP} 
(for other $u$, the assumption is trivially satisfied),
reads as
\begin{equation}\label{eq:new_CQ_NLP}
	\left.
		\begin{aligned}
			\sum\nolimits_{i=1}^m\lambda_i\nabla^2g_i(\bar x)u
				+\nabla g(\bar x)^*\tilde\lambda=0,\,\nabla g(\bar x)^*\lambda=0,\\
			\forall i\in I^0(\bar x,u)\colon\,\lambda_i\geq 0,\\
			\forall i\in I^{00}(\bar x,u,\lambda)\colon\,\tilde\lambda_i\geq 0,\\
			\forall i\in I^-(\bar x,u)\cup I^-(\bar x)\colon\,\lambda_i=\tilde\lambda_i=0
		\end{aligned}
	\right\}
	\quad
	\Longrightarrow
	\quad
	\lambda=0
\end{equation}
where we used
\begin{align*}
	I^0(\bar x,u)&:=\{i\in I(\bar x)\,|\,\nabla g_i(\bar x) u=0\},\\
	I^-(\bar x,u)&:=\{i\in I(\bar x)\,|\,\nabla g_i(\bar x) u<0\},\\
	I^{00}(\bar x,u,\lambda)&:=\{i\in I^0(\bar x,u)\,|\,\lambda_i=0\}.
\end{align*}
Now, we can distill from \cref{Pro:M-stat_via_second_order}\,\ref{item:pseudo_stationarity_constraint_maps_polyhedral}
that whenever $\bar x$ is a local minimizer of the associated problem \eqref{eq:nonsmooth_problem},
and, for each $u\in\mathbb S_{\mathbb X}$ satisfying \eqref{eq:linearization_cone_NLP}, we have \eqref{eq:new_CQ_NLP},
then either $\bar x$ is M-stationary, i.e., we find $\lambda\in\R^m$ such that
\[
	\begin{aligned}
		&0\in\partial\varphi(\bar x)+\nabla g(\bar x)^*\lambda,\\
		&\forall i\in I(\bar x)\colon\,\lambda_i\geq 0,\\
		&\forall i\in I^-(\bar x)\colon\,\lambda_i=0,
	\end{aligned}
\]
or we find $u\in\mathbb S_{\mathbb X}$ satisfying \eqref{eq:linearization_cone_NLP} and $\lambda,\tilde\lambda\in\R^m$ such that
\begin{equation}\label{eq:mixed_M_stat_NLP}
	\begin{aligned}
		&0\in\partial\varphi(\bar x)+\sum\nolimits_{i=1}^m\lambda_i\nabla^2 g_i(\bar x)u+\nabla g(\bar x)^*\tilde\lambda,\,\nabla g(\bar x)^*\lambda = 0,\\
		&\forall i\in I^0(\bar x,u)\colon\,\lambda_i\geq 0,\\
		&\forall i\in I^{00}(\bar x,u,\lambda)\colon\,\tilde\lambda_i\geq 0,\\
		&\forall i\in I^-(\bar x,u)\cup I^-(\bar x)\colon\,\lambda_i=\tilde\lambda_i=0.
	\end{aligned}
\end{equation}
Keeping \cref{rem:upper_estimate_pseudo_coderivative_constraint_maps} in mind, this is
the same stationarity condition which could be obtained from
\cref{Pro:M-stat_via_second_order}\,\ref{item:pseudo_stationarity_constraint_maps}.

\begin{example}\label{ex:GCQ_violated}
	We consider
	\[
		\min\{x\,|\,x^2\leq 0\}
	\]
	and its uniquely determined feasible point $\bar x:=0$ 
	where even Guignard's constraint qualification is violated, see \cite{GouldTolle1971}.
	Particularly, the M-stationarity conditions 
	which equal the KKT conditions in this situation, do not hold.
	In this situation, however, $u=\pm 1$ satisfies \eqref{eq:linearization_cone_NLP} and due to
	\[
		2\lambda u=0,\,\lambda\geq 0\quad\Longrightarrow\quad\lambda =0,
	\]
	\eqref{eq:new_CQ_NLP} holds. 
	The stationarity conditions \eqref{eq:mixed_M_stat_NLP} are satisfied, 
	e.g., with $u:=-1$, $\lambda:=1/2$, and $\tilde\lambda:=0$.
\end{example}

Now, let us focus on \cref{Pro:M-stat_via_second_order}\,\ref{item:pseudo_stationarity_constraint_maps_polyhedral_refined}.
For $u\in\mathbb S_{\mathbb X}$ satisfying \eqref{eq:linearization_cone_NLP}, we find
\[
	\mathbf T(u)
	=
	\left\{
		y\in\R^m\,\middle|\,
		\begin{aligned}
			&\forall i\in I^0(\bar x,u)\colon\,y_i\leq 0\\
			&\forall i\in\{\ell+1,\ldots,m\}\colon\,y_i=0
		\end{aligned}
	\right\}.
\]
For $s\in\mathbb X$, $v\in\R^m$, and $\lambda\in\R^m$, we further define
\begin{align*}
	I^0_0(\bar x,u,s,v)
	&:=
	\{i\in I^0(\bar x,u)\,|\,\nabla g_i(\bar x) s+\tfrac12\nabla^2g_i(\bar x)[u,u]=v_i\},\\
	I^0_-(\bar x,u,s,v)
	&:=
	\{i\in I^0(\bar x,u)\,|\,\nabla g_i(\bar x) s+\tfrac12 \nabla^2g_i(\bar x)[u,u]<v_i\},\\
	I^{00}_0(\bar x,u,s,v,\lambda)
	&:=
	\{i\in I^0_0(\bar x,u,s,v)\,|\,\lambda_i=0\}.
\end{align*}
Then \eqref{eq:CQ_pseudo_subregularity_polyhedral_II}
demands that for each $u\in\mathbb S_{\mathbb X}$ satisfying \eqref{eq:linearization_cone_NLP}, we find
\begin{equation}\label{eq:new_CQ_NLP_refined}
	\left.
		\begin{aligned}
		&\sum\nolimits_{i=1}^m\lambda_i\nabla^2g_i(\bar x)u
			+\nabla g(\bar x)^*\tilde\lambda=0,\,\nabla g(\bar x)^*\lambda=0,\\
		&\forall i\in I^0(\bar x,u)\colon\,
			\nabla g_i(\bar x) s+\tfrac12\nabla^2g_i(\bar x)[u,u]\leq 0,\\
		&\forall i\in\{\ell+1,\ldots,m\}\colon\,
			\nabla g_i(\bar x) s+\tfrac12\nabla^2g_i(\bar x)[u,u]= 0,\\
		&\forall i\in I^0_0(\bar x,u,s,0)\colon\,\lambda_i\geq 0,\\
		&\forall i\in I^{00}_0(\bar x,u,s,0,\lambda)\colon\,\tilde\lambda_i\geq 0,\\
		&\forall i\in I^0_-(\bar x,u,s,0)\cup I^-(\bar x,u)\cup I^-(\bar x)\colon\,\lambda_i=\tilde\lambda_i=0
		\end{aligned}
	\right\}
	\quad
	\Longrightarrow
	\quad
	\lambda=0.
\end{equation}
Whenever $\bar x$ is a local minimizer of the associated problem \eqref{eq:nonsmooth_problem}, and, for each $u\in\mathbb S_{\mathbb X}$ satisfying
\eqref{eq:linearization_cone_NLP}, we have \eqref{eq:new_CQ_NLP_refined}, then either $\bar x$ is M-stationary or we find $u\in\mathbb S_{\mathbb X}$
satisfying \eqref{eq:linearization_cone_NLP} and $s\in\mathbb X$ as well as $\lambda,\tilde\lambda,v\in\R^m$ such that
\begin{equation}\label{eq:mixed_M_stat_NLP_refined}
	\begin{aligned}
		&0\in\partial\varphi(\bar x)+\sum\nolimits_{i=1}^m\lambda_i\nabla^2g_i(\bar x)u
			+\nabla g(\bar x)^*\tilde\lambda,\,
			\nabla g(\bar x)^*\lambda=0,\,\innerprod{ v }{\lambda}\geq 0,\\
		&\forall i\in I^0(\bar x,u)\colon\,
			\nabla g_i(\bar x) s+\tfrac12\nabla^2g_i(\bar x)[u,u]\leq v,\\
		&\forall i\in\{\ell+1,\ldots,m\}\colon\,
			\nabla g_i(\bar x) s+\tfrac12\nabla^2g_i(\bar x)[u,u]= v,\\
		&\forall i\in I^0_0(\bar x,u,s,v)\colon\,\lambda_i\geq 0,\\
		&\forall i\in I^{00}_0(\bar x,u,s,v,\lambda)\colon\,\tilde\lambda_i\geq 0,\\
		&\forall i\in I^0_-(\bar x,u,s,v)\cup I^-(\bar x,u)\cup I^-(\bar x)\colon\,
			\lambda_i=\tilde\lambda_i=0.
	\end{aligned}
\end{equation}
By definition of the appearing index sets, we directly see that \eqref{eq:new_CQ_NLP_refined} is not stronger than \eqref{eq:new_CQ_NLP}
while validity of the stationarity condition \eqref{eq:mixed_M_stat_NLP_refined} implies validity of \eqref{eq:mixed_M_stat_NLP}.

We note that these necessary optimality conditions are different from those ones from
\cite{Avakov1989,AvakovArutunovIzmailov2007,Gfrerer2007,Gfrerer2014a}
which hold without further assumptions but comprise a (potentially vanishing)
leading multiplier associated with the objective function.

\subsection{Mathematical programs with complementarity constraints}\label{sec:MPCCs}

Let us introduce
\[
	\mathcal C:=(\R_+\times\{0\})\cup(\{0\}\times\R_+),
\]
the so-called complementarity angle.
For twice continuously differentiable data functions $G,H\colon\mathbb X\to\R^m$ 
with components $G_1,\ldots,G_m\colon\mathbb X\to\R$ and $H_1,\ldots,H_m\colon\mathbb X\to\R$, 
we address the constraint region given by
\begin{equation}\label{eq:MPCC}\tag{MPCC}
	(G_i(x),H_i(x))\in\mathcal C\quad i\in I
\end{equation}
where $I:=\{1,\ldots,m\}$.
The latter is distinctive for so called 
\emph{mathematical programs with complementarity constraints}
which have been studied intensively throughout the last decades,
see e.g.\ \cite{LuoPangRalph1996,OutrataKocvaraZowe1998} for some classical references.
Again, we observe that \eqref{eq:MPCC} can be formulated via a constraint map using 
$D:=\mathcal C^m$.
Note that standard inequality and equality constraints 
can be added without any difficulties
due to \cref{lem:product_rule_tangents_polyhedral_sets,lem:product_rule_graphical_derivative} by
applying our findings from \cref{sec:NLPs}. Here, we omit them for brevity of presentation. 

Fix some feasible point $\bar x\in\mathbb X$ of \eqref{eq:MPCC}.
A critical direction $u\in\mathbb S_{\mathbb X}$ 
of the associated problem \eqref{eq:nonsmooth_problem} necessarily needs to satisfy
\begin{equation}\label{eq:linearization_cone_MPCC}
	\begin{aligned}
		\nabla G_i(\bar x) u&=0\quad& &i\in I^{0+}(\bar x),&\\
		\nabla H_i(\bar x) u&=0\quad&	&i\in I^{+0}(\bar x),&\\
		(\nabla G_i(\bar x) u,\nabla H_i(\bar x) u)&\in\mathcal C&&i\in I^{00}(\bar x),&
	\end{aligned}
\end{equation}
where we used the well-known index sets
\begin{align*}
	I^{0+}(\bar x)&:=\{i\in I\,|\,G_i(\bar x)=0,\,H_i(\bar x)>0\},\\
	I^{+0}(\bar x)&:=\{i\in I\,|\,G_i(\bar x)>0,\,H_i(\bar x)=0\},\\
	I^{00}(\bar x)&:=\{i\in I\,|\,G_i(\bar x)=0,\,H_i(\bar x)=0\}.
\end{align*}
Noting that \eqref{eq:CQ_pseudo_subregularity_Ia} and \eqref{eq:CQ_pseudo_subregularity_Ib}
reduce to \eqref{eq:CQ_pseudo_subregularity_I_polyhedral} by 
\cref{rem:upper_estimate_pseudo_coderivative_constraint_maps}\,\ref{item:graphical_subderivative_assumption_polyhedral},
with the aid of \cref{lem:some_properties_of_polyhedral_sets,lem:product_rule_tangents_polyhedral_sets,lem:product_rule_graphical_derivative},
we now can characterize 
\eqref{eq:CQ_pseudo_subregularity_I_polyhedral} and \eqref{eq:CQ_pseudo_subregularity_II}.
Based on the representation
\[
		\gph\mathcal N_{\mathcal C}
		=
		(\R_+\times\{0\}\times\{0\}\times\R)\cup(\{0\}\times\R_+\times\R\times\{0\})\cup(\{0\}\times\{0\}\times\R_-\times\R_-),
	\]
some elementary calculations show
\begin{align*}
		D\mathcal N_{\mathcal C}((a,b),(\mu,\nu))(v)
		=
		\begin{cases}
			\{0\}\times\R	&	a>0,\,b=\mu=0,\,v_2=0,\\
			\R\times\{0\}	&	a=\nu=0,\,b>0,\,v_1=0,\\
			\R^2			&	a=b=0,\,\mu,\nu<0,\,v=0,\\
			\{0\}\times\R	&	a=b=\mu=0,\,\nu<0,\,v_1>0,\,v_2=0,\\
			\R_-\times\R	&	a=b=\mu=0,\,\nu<0,\,v=0,\\
			\R\times\{0\}	&	a=b=\nu=0,\,\mu<0,\,v_1=0,\,v_2>0,\\
			\R\times\R_-	&	a=b=\nu=0,\,\mu<0,\,v=0,\\
			\{0\}\times\R	&	a=b=\mu=0,\,\nu>0,\,v_1\geq 0,\,v_2=0,\\
			\R\times\{0\}	&	a=b=\nu=0,\,\mu>0,\,v_1=0,\,v_2\geq 0,\\
			\{0\}\times\R	&	a=b=\mu=\nu=0,\,v_1>0,\,v_2=0,\\
			\R\times\{0\}	&	a=b=\mu=\nu=0,\,v_1=0,\,v_2>0,\\
			\mathcal N_{\mathcal C}(0)	&	a=b=\mu=\nu=0,\,v=0,\\
			\varnothing		&	\text{otherwise}
		\end{cases}
\end{align*}
for arbitrary $((a,b),(\mu,\nu))\in\gph\mathcal N_{\mathcal C}$ and $v\in\R^2$.
Consequently, for $u\in\mathbb S_{\mathbb X}$ 
satisfying \eqref{eq:linearization_cone_MPCC}, we need to ensure
\[
	\left.
		\begin{aligned}
		&\nabla G(\bar x)^*\mu+\nabla H(\bar x)^*\nu=0,\,
		\nabla G(\bar x)^*\tilde\mu+\nabla H(\bar x)^*\tilde\nu=0,\\
		&\forall i\in I^{+0}(\bar x)\cup I^{00}_{+0}(\bar x,u)\colon\,\mu_i=0,\\
		&\forall i\in I^{0+}(\bar x)\cup I^{00}_{0+}(\bar x,u)\colon\,\nu_i=0,\\
		&\forall i\in I^{00}_{00}(\bar x,u)\colon\,\mu_i,\nu_i\leq 0\,\lor\,\mu_i\nu_i=0,\\
		&\forall i\in I\colon\,(\tilde\mu_i,\tilde\nu_i)
		\in D\mathcal N_{\mathcal C}
		((\bar G_i,\bar H_i),(\mu_i,\nu_i))(\nabla\bar G_iu,\nabla \bar H_iu)
		\end{aligned}
	\right\}
	\quad
	\Longrightarrow
	\quad
	\tilde\mu=\tilde\nu=0
\]
and
\begin{equation}\label{eq:CQ_for_mixed_order_stationarity_MPCC}
	\left.
		\begin{aligned}
		&\nabla G(\bar x)^*\mu+\nabla H(\bar x)^*\nu=0,\\
		&\sum\nolimits_{i=1}^m\bigl(\mu_i\nabla^2G_i(\bar x)+\nu_i\nabla^2H_i(\bar x)\bigr)u
		\\
		&\qquad\qquad
		+\nabla G(\bar x)^*\tilde\mu+\nabla H(\bar x)^*\tilde\nu=0,\\
		&\forall i\in I^{+0}(\bar x)\cup I^{00}_{+0}(\bar x,u)\colon\,\mu_i=0,\\
		&\forall i\in I^{0+}(\bar x)\cup I^{00}_{0+}(\bar x,u)\colon\,\nu_i=0,\\
		&\forall i\in I^{00}_{00}(\bar x,u)\colon\,\mu_i,\nu_i\leq 0\,\lor\,\mu_i\nu_i=0,\\
		&\forall i\in I\colon\,(\tilde\mu_i,\tilde\nu_i)
		\in D\mathcal N_{\mathcal C}
		((\bar G_i,\bar H_i),(\mu_i,\nu_i))(\nabla\bar G_iu,\nabla \bar H_iu)
		\end{aligned}
	\right\}
	\quad
	\Longrightarrow
	\quad
	\mu=\nu=0
\end{equation}
where, for each $i\in I$, we used $\bar G_i:=G_i(\bar x)$, $\bar H_i:=H_i(\bar x)$, 
$\nabla\bar G_iu:=\nabla G_i(\bar x) u$, and
$\nabla\bar H_iu:=\nabla H_i(\bar x) u$ for brevity as well as the index sets
\begin{align*}
	I^{00}_{+0}(\bar x,u)&:=\{i\in I^{00}(\bar x)\,|\,\nabla \bar G_iu>0,\,\nabla \bar H_iu=0\},\\
	I^{00}_{0+}(\bar x,u)&:=\{i\in I^{00}(\bar x)\,|\,\nabla \bar G_iu=0,\,\nabla \bar H_iu>0\},\\
	I^{00}_{00}(\bar x,u)&:=\{i\in I^{00}(\bar x)\,|\,\nabla \bar G_iu=0,\,\nabla \bar H_iu=0\}.
\end{align*}
The assertion of \cref{Pro:M-stat_via_second_order}\,\ref{item:pseudo_stationarity_constraint_maps} now yields that whenever 
$\bar x$ is a local minimizer for the associated problem \eqref{eq:nonsmooth_problem} and for each
$u\in\mathbb S_{\mathbb X}$ satisfying \eqref{eq:linearization_cone_MPCC}, the two above implications hold, then $\bar x$ is either
M-stationary, i.e., there are multipliers $\mu,\nu\in\R^m$ satisfying
\begin{align*}
	&0\in\partial \varphi(\bar x)+\nabla G(\bar x)^*\mu+\nabla H(\bar x)^*\nu,\\
	&\forall i\in I^{+0}(\bar x)\colon\,\mu_i=0,\\
	&\forall i\in I^{0+}(\bar x)\colon\,\nu_i=0,\\
	&\forall i\in I^{00}(\bar x)\colon\,\mu_i,\nu_i\leq 0\,\lor\,\mu_i\nu_i=0,
\end{align*}
or we find $u\in\mathbb S_{\mathbb X}$ satisfying \eqref{eq:linearization_cone_MPCC} as well as multipliers $\mu,\nu,\tilde\mu,\tilde\nu\in\R^m$ satisfying
\begin{equation}\label{eq:mixed_order_stationarity_MPCC}
	\begin{aligned}
	&0\in\partial \varphi(\bar x)
		+\sum\nolimits_{i=1}^m\bigl(\mu_i\nabla^2G_i(\bar x)+\nu_i\nabla^2H_i(\bar x)\bigr)u
		+\nabla G(\bar x)^*\tilde\mu+\nabla H(\bar x)^*\tilde\nu,\\
	&0=\nabla G(\bar x)^*\mu+\nabla H(\bar x)^*\nu,\\
	&\forall i\in I^{+0}(\bar x)\cup I^{00}_{+0}(\bar x,u)\colon\,\mu_i=0,\\
	&\forall i\in I^{0+}(\bar x)\cup I^{00}_{0+}(\bar x,u)\colon\,\nu_i=0,\\
	&\forall i\in I^{00}_{00}(\bar x,u)\colon\,\mu_i,\nu_i\leq 0\,\lor\,\mu_i\nu_i=0,\\
	&\forall i\in I\colon\,(\tilde\mu_i,\tilde\nu_i)
	\in D\mathcal N_{\mathcal C}
	((\bar G_i,\bar H_i),(\mu_i,\nu_i))(\nabla\bar G_iu,\nabla \bar H_iu).
	\end{aligned}
\end{equation}

In order to characterize \eqref{eq:CQ_pseudo_subregularity_polyhedral_I},
we observe that
\begin{align*}
	\mathcal N_{\mathcal T_{\mathcal C}(\bar G_i,\bar H_i)}(\nabla\bar G_iu,\nabla\bar H_iu)
	=
	\begin{cases}
		\{0\}\times\R	&i\in I^{+0}(\bar x)\cup I^{00}_{+0}(\bar x,u),\\
		\R\times\{0\}	&i\in I^{0+}(\bar x)\cup I^{00}_{0+}(\bar x,u),\\
		\mathcal N_{\mathcal C}(0)	&i\in I^{00}_{00}(\bar x,u)
	\end{cases}
\end{align*}
is valid for each $i\in I$. 
For each pair $(\mu_i,\nu_i)\in\mathcal N_{\mathcal T_{\mathcal C}(\bar G_i,\bar H_i)}(\nabla\bar G_iu,\nabla\bar H_iu)$,
elementary calculations show
\begin{align*}
	\mathcal T_{\mathcal N_{\mathcal T_{\mathcal C}(\bar G_i,\bar H_i)}(\nabla\bar G_iu,\nabla\bar H_iu)}(\mu_i,\nu_i)
	&=
	\begin{cases}
		\{0\}\times\R	&i\in I^{+0}(\bar x)\cup I^{00}_{+0}(\bar x,u),\\
		\R\times\{0\}	&i\in I^{0+}(\bar x)\cup I^{00}_{0+}(\bar x,u),\\
		\R^2			&i\in I^{00}_{00}(\bar x,u),\,\mu_i<0,\,\nu_i<0,\\
		\R_-\times\R	&i\in I^{00}_{00}(\bar x,u),\,\mu_i=0,\,\nu_i<0,\\
		\R\times\R_-	&i\in I^{00}_{00}(\bar x,u),\,\mu_i<0,\,\nu_i=0,\\
		\{0\}\times\R	&i\in I^{00}_{00}(\bar x,u),\,\mu_i=0,\,\nu_i>0,\\
		\R\times\{0\}	&i\in I^{00}_{00}(\bar x,u),\,\mu_i>0,\,\nu_i=0,\\
		\mathcal N_{\mathcal C}(0)	&i\in I^{00}_{00}(\bar x,u),\,\mu_i=\nu_i=0
	\end{cases}
	\\
	&=
	D\mathcal N_{\mathcal C}((\bar G_i,\bar H_i),(\mu_i,\nu_i))(\nabla \bar G_i,\nabla \bar H_i).
\end{align*}
Thus, due to \cref{Pro:M-stat_via_second_order}\,\ref{item:pseudo_stationarity_constraint_maps_polyhedral}, 
validity of \eqref{eq:CQ_for_mixed_order_stationarity_MPCC}
for each $u\in\mathbb S_{\mathbb X}$ satisfying \eqref{eq:linearization_cone_MPCC} 
is already enough to infer that whenever $\bar x$ is a local minimizer,
then it is either M-stationary or there are $u\in\mathbb S_{\mathbb X}$ satisfying \eqref{eq:linearization_cone_MPCC} 
and multipliers $\mu,\nu,\tilde\mu,\tilde\nu\in\R^m$ 
solving the stationarity conditions \eqref{eq:mixed_order_stationarity_MPCC}.

For brevity, we abstain from presenting the results from 
\cref{Pro:M-stat_via_second_order}\,\ref{item:pseudo_stationarity_constraint_maps_polyhedral_refined}
for the constraint system \eqref{eq:MPCC}.

\subsection{Semidefinite programming}

Let us consider the Hilbert space $\SSS_m$ of all real symmetric matrices equipped with the
standard (Frobenius) inner product. We denote by $\SSS_m^+$ and $\SSS_m^-$ the cone of all
positive and negative semidefinite matrices in $\SSS_m$, respectively.
The foundations of variational analysis in this space can be found, 
e.g., in \cite[Section~5.3]{BonnansShapiro2000}.
For some twice continuously differentiable mapping $g\colon\mathbb X\to\SSS_m$, 
we investigate the
constraint system
\begin{equation}\label{eq:semidefinite_NLP}\tag{NLSD}
	g(x)\in\SSS_m^+.
\end{equation} 
It is well known that the closed, convex cone $\SSS_m^+$ is not polyhedral.
Nevertheless, the constraint \eqref{eq:semidefinite_NLP} can be encoded 
via a constraint map.

Let $\bar x\in\mathbb X$ be feasible to \eqref{eq:semidefinite_NLP} and, 
for some $u\in\mathbb S_{\mathbb X}$, fix
$\Omega\in\mathcal N_{\SSS_m^+}(g(\bar x);\nabla g(\bar x)u)$. 
For later use, fix an orthogonal matrix 
$\mathbf P\in\R^{m\times m}$ and a diagonal matrix $\mathbf\Lambda\in\R^{m\times m}$
whose diagonal elements $\lambda_1,\ldots,\lambda_m$ are ordered non-increasingly 
such that $g(\bar x)+\Omega=\mathbf P\mathbf\Lambda\mathbf P^\top$. The index sets
corresponding to the positive, zero, and negative entries on the main diagonal 
of $\mathbf\Lambda$
are denoted by $\alpha$, $\beta$, and $\gamma$, respectively.
We emphasize that, here and throughout the subsection,
$\alpha$ is a constant index set while
$\beta$ and $\gamma$ depend on the precise choice of $\Omega$.
Subsequently, we use the notation 
$\mathbf Q^{\mathbf P}:=\mathbf P^\top \mathbf Q\mathbf P$ and 
$\mathbf Q^{\mathbf P}_{IJ}:=(\mathbf Q^{\mathbf P})_{IJ}$
for each matrix $\mathbf Q\in\mathcal S_m$ 
and index sets $I,J\subset\{1,\ldots,m\}$
where $\mathbf M_{IJ}$ is the submatrix of $\mathbf M\in\mathcal S_m$
which possesses only those rows and columns of $\mathbf M$ whose
indices can be found in $I$ and $J$, respectively. 

The above constructions yield
\[
	g(\bar x)=\mathbf P\max(\mathbf\Lambda,\mathbf O)\mathbf P^\top,\qquad
	\Omega=\mathbf P\min(\mathbf\Lambda,\mathbf O)\mathbf P^\top
\]
where $\max$ and $\min$ have to be understood in entrywise fashion and $\OOO$ is an all-zero matrix of appropriate dimensions.
Due to
\[
	\nabla g(\bar x)u\in\mathcal T_{\mathcal S_m^+}(g(\bar x))
	=
	\left\{
		\QQQ\in\mathcal S_m\,\middle|\,
		\QQQ^{\PPP}_{\beta\cup\gamma,\beta\cup\gamma}\in\mathcal S_{|\beta\cup\gamma|}^+
	\right\},
\]
we find
\begin{align*}
	0
	&=
	\innerprod{\Omega}{\nabla g(\bar x)u}
	=
	\trace(\Omega\,\nabla g(\bar x)u)
	=
	\trace(\PPP\min(\LLL,\mathbf O)\PPP^\top \PPP[\nabla g(\bar x)u]^{\PPP}\PPP^\top)
	\\
	&=
	\trace(\min(\LLL,\mathbf O)[\nabla g(\bar x)u]^{\PPP})
	=
	\sum_{i\in\gamma}\underbrace{\lambda_i}_{<0}\,\underbrace{[\nabla g(\bar x)u]^{\PPP}_{i,i}}_{\geq 0}
\end{align*}
which directly gives us $[\nabla g(\bar x)u]^\PPP_{\beta\gamma}=\OOO$, $[\nabla g(\bar x)u]^\PPP_{\gamma\gamma}=\OOO$, and $[\nabla g(\bar x)u]^\PPP_{\beta\beta}\in\SSS_{|\beta|}^+$.
Furthermore, we note
\[
	\mathcal N_{\SSS^+_m}(g(\bar x))
	=
	\left\{
		\tilde\Omega\in\SSS_m\,\middle|\,
		\tilde\Omega^\PPP_{\alpha\alpha}=\OOO,\,\tilde\Omega^\PPP_{\alpha\beta}=\OOO,\,\tilde\Omega^{\PPP}_{\alpha\gamma}=\OOO,\,
		\tilde\Omega^\PPP_{\beta\cup\gamma,\beta\cup\gamma}\in\SSS^-_{|\beta\cup\gamma|}
	\right\}.
\]
Finally, let $\Xi_{\alpha\gamma}\in\R^{|\alpha|\times|\gamma|}$ be the matrix given by
\[
	\forall i\in \alpha\,\forall j\in\gamma\colon\quad
	[\Xi_{\alpha\gamma}]_{ij}:=-\frac{\lambda_j}{\lambda_i}.
\]
It is a well-known result that the projection onto $\SSS_m^+$ is directionally
differentiable.
With the aid of \cref{lem:graphical_derivatives_of_normal_cone_map} and \cite[Corollary~3.1]{WuZhangZhang2014}, we find
\begin{align*}
	D\mathcal N_{\SSS_m^+}(g(\bar x),\Omega)(\nabla g(\bar x)u)
	&=
	\left\{\tilde\Omega\in\SSS_m\,\middle|\,
		\begin{aligned}
			&\tilde\Omega^\PPP_{\alpha\alpha}=\OOO,\,\tilde\Omega^\PPP_{\alpha\beta}=\OOO,\,
			\tilde\Omega^\PPP_{\alpha\gamma}=\Xi_{\alpha\gamma}\bullet[\nabla g(\bar x)u]^\PPP_{\alpha\gamma},\\
			&\tilde\Omega^\PPP_{\beta\beta}\in\SSS_{|\beta|}^-,\,\innerprod{\tilde\Omega^\PPP_{\beta\beta}}{[\nabla g(\bar x)u]^\PPP_{\beta\beta}}=0
		\end{aligned}
	\right\},
\end{align*}
and if $\nabla g(\bar x)u\neq\OOO$, we find
\begin{align*}
	D_\textup{sub}\mathcal N_{\SSS_m^+}(g(\bar x),\Omega)\left(\frac{\nabla g(\bar x)u}{\norm{\nabla g(\bar x)u}}\right)
	&\subset
	\left\{\tilde\Omega\in\SSS_m\,\middle|\,
		\begin{aligned}
			&\tilde\Omega^\PPP_{\alpha\alpha}=\OOO,\,\tilde\Omega^\PPP_{\alpha\beta}=\OOO,\,
			\tilde\Omega^\PPP_{\alpha\gamma}=\OOO,\\
			&\tilde\Omega^\PPP_{\beta\beta}\in\SSS_{|\beta|}^-,\,
			\innerprod{\tilde\Omega^\PPP_{\beta\beta}}{[\nabla g(\bar x)u]^\PPP_{\beta\beta}}=0
		\end{aligned}
	\right\}.
\end{align*}
Above, $\bullet$ represents the \emph{Hadamard}, i.e., entrywise product.
Note that validity of the final orthogonality condition in the estimate 
for the graphical subderivative follows from
\cref{lem:graphical_derivatives_of_normal_cone_map} since 
$\tilde\Omega\in D_\textup{sub}\mathcal N_{\SSS_m^+}(g(\bar x),\Omega)(\nabla g(\bar x)u/\norm{\nabla g(\bar x)u})$ 
and $\norm{\nabla g(\bar x)u}>0$ yield
\begin{align*}
	0
	&\leq
	\innerprod{\tilde\Omega}{\nabla g(\bar x)u}
	=
	\trace\bigl(\tilde\Omega\,\nabla g(\bar x)u\bigr)
	=
	\trace\bigl(\PPP\tilde\Omega^{\PPP}\PPP^\top\PPP[\nabla g(\bar x)u]^\PPP\PPP^\top\bigr)
	\\
	&=
	\trace\bigl(\tilde\Omega^\PPP[\nabla g(\bar x)u]^\PPP\bigr)
	=
	\trace\bigl(\tilde\Omega^\PPP_{\beta\beta}[\nabla g(\bar x)u]^\PPP_{\beta\beta}\bigr)
	\leq 
	0
\end{align*}
due to $\tilde\Omega^\PPP_{\alpha\alpha}=\OOO$, $\tilde\Omega^\PPP_{\alpha\beta}=\OOO$, $\tilde\Omega^\PPP_{\alpha\gamma}=\OOO$,
$\tilde\Omega^\PPP_{\beta\beta}\in\SSS_{|\beta|}^-$, $[\nabla g(\bar x)u]^\PPP_{\beta\gamma}=\OOO$, $[\nabla g(\bar x)u]^\PPP_{\gamma\gamma}=\OOO$,
and $[\nabla g(\bar x)u]^\PPP_{\beta\beta}\in\SSS_{|\beta|}^+$.
Thus, for each $u\in\mathbb S_{\mathbb X}$, \eqref{eq:CQ_pseudo_subregularity_II} takes the form
\[
	\left.
		\begin{aligned}
			\nabla g(\bar x)^*\Omega=0,\,
			\nabla^2\innerprod{\Omega}{g}(\bar x)(u)+\nabla g(\bar x)^*\tilde\Omega=0,\\
			\Omega^\PPP_{\alpha\alpha}=\OOO,\,\Omega^\PPP_{\alpha\beta}=\OOO,\,
			\Omega^{\PPP}_{\alpha\gamma}=\OOO,\,
			\Omega^\PPP_{\beta\cup\gamma,\beta\cup\gamma}\in\SSS_{|\beta\cup\gamma|}^-,\\
			[\nabla g(\bar x)u]^\PPP_{\beta\gamma}=\OOO,\,
			[\nabla g(\bar x)u]^\PPP_{\gamma\gamma}=\OOO,\,
			[\nabla g(\bar x)u]^\PPP_{\beta\beta}\in\SSS_{|\beta|}^+,\\
			\tilde\Omega^\PPP_{\alpha\alpha}=\OOO,\,\tilde\Omega^\PPP_{\alpha\beta}=\OOO,\,
			\tilde\Omega^\PPP_{\alpha\gamma}
			=
			\Xi_{\alpha\gamma}\bullet[\nabla g(\bar x)u]^\PPP_{\alpha\gamma},\\
			\tilde\Omega^\PPP_{\beta\beta}\in\SSS_{|\beta|}^-,\,
			\innerprod{\tilde\Omega^\PPP_{\beta\beta}}{[\nabla g(\bar x)u]^\PPP_{\beta\beta}}=0
		\end{aligned}
	\right\}
	\quad
	\Longrightarrow\quad
	\Omega=\OOO,
\]
while \eqref{eq:CQ_pseudo_subregularity_Ia} and \eqref{eq:CQ_pseudo_subregularity_Ib} 
(the latter in case $\nabla g(\bar x)u\neq\OOO$) are both implied by
\[
	\left.
		\begin{aligned}
			\nabla g(\bar x)^*\Omega=0,\,\nabla g(\bar x)^*\tilde\Omega=0,\\
			\Omega^\PPP_{\alpha\alpha}=\OOO,\,\Omega^\PPP_{\alpha\beta}=\OOO,\,
			\Omega^{\PPP}_{\alpha\gamma}=\OOO,\,
			\Omega^\PPP_{\beta\cup\gamma,\beta\cup\gamma}\in\SSS_{|\beta\cup\gamma|}^-,\\
			[\nabla g(\bar x)u]^\PPP_{\beta\gamma}=\OOO,\,
			[\nabla g(\bar x)u]^\PPP_{\gamma\gamma}=\OOO,\,
			[\nabla g(\bar x)u]^\PPP_{\beta\beta}\in\SSS_{|\beta|}^+,\\
			\tilde\Omega^\PPP_{\alpha\alpha}=\OOO,\,\tilde\Omega^\PPP_{\alpha\beta}=\OOO,\,
			\tilde\Omega^\PPP_{\alpha\gamma}=\OOO,\\
			\tilde\Omega^\PPP_{\beta\beta}\in\SSS_{|\beta|}^-,\,
			\innerprod{\tilde\Omega^\PPP_{\beta\beta}}{[\nabla g(\bar x)u]^\PPP_{\beta\beta}}=0
		\end{aligned}
	\right\}
	\quad
	\Longrightarrow\quad
	\tilde\Omega=\OOO.
\]
In case where $\bar x$ is a local minimizer of the associated problem \eqref{eq:nonsmooth_problem}, 
validity of these conditions for each $u\in\mathbb S_{\mathbb X}$ guarantees that
$\bar x$ is either M-stationary (we omit stating this well-known system here) 
or we find $u\in\mathbb S_{\mathbb X}$ 
and $\Omega,\tilde\Omega\in\SSS_m$ such that
\begin{align*}
			&0\in\partial \varphi(\bar x)+\nabla^2\innerprod{\Omega}{g}(\bar x)(u)
			+\nabla g(\bar x)^*\tilde\Omega,\,\nabla g(\bar x)^*\Omega=0,\\
			&\Omega^\PPP_{\alpha\alpha}=\OOO,\,\Omega^\PPP_{\alpha\beta}=\OOO,\,
			\Omega^{\PPP}_{\alpha\gamma}=\OOO,\,
			\Omega^\PPP_{\beta\cup\gamma,\beta\cup\gamma}\in\SSS_{|\beta\cup\gamma|}^-,\\
			&[\nabla g(\bar x)u]^\PPP_{\beta\gamma}=\OOO,\,
			[\nabla g(\bar x)u]^\PPP_{\gamma\gamma}=\OOO,\,
			[\nabla g(\bar x)u]^\PPP_{\beta\beta}\in\SSS_{|\beta|}^+,\\
			&\tilde\Omega^\PPP_{\alpha\alpha}=\OOO,\,
			\tilde\Omega^\PPP_{\alpha\beta}=\OOO,\,
			\tilde\Omega^\PPP_{\alpha\gamma}
			=\Xi_{\alpha\gamma}\bullet[\nabla g(\bar x)u]^\PPP_{\alpha\gamma},\\
			&\tilde\Omega^\PPP_{\beta\beta}\in\SSS_{|\beta|}^+,\,
			\innerprod{\tilde\Omega^\PPP_{\beta\beta}}{[\nabla g(\bar x)u]^\PPP_{\beta\beta}}=0.
\end{align*}
	
\section{Concluding remarks}\label{sec:conclusions}

In this paper, we enriched the general concept of approximate stationarity with the aid
of tools from directional limiting variational analysis. Based on our main result
\cref{thm:higher_order_directional_asymptotic_stationarity}, we were in position to
obtain new mixed-order necessary optimality conditions and associated constraint
qualifications which characterize local minimizers. Some new upper estimates
for the second-order directional pseudo-coderivative of constraint mappings were
successfully employed to make these results fully explicit in the presence of
geometric constraints. Our findings also gave rise to new applicable qualification conditions
guaranteeing M-stationarity of local minimizers,
and in Part B of this paper, we show that these conditions 
are not stronger than FOSCMS and the Second-Order Sufficient Condition for Metric Subregularity. 
It is a topic of future research to specify these results for other constraint regions,
e.g., given by equilibrium conditions.
Similarly, it seems desirable to further develop the calculus for pseudo-coderivatives.
Based on \cref{thm:directional_asymptotic_stationarity} it is also possible to introduce
new mild \emph{sequential} constraint qualifications for M-stationarity which,
in contrast to the concept from \cite{Mehlitz2020a}, depend on critical directions.
This is also done in Part B of this paper.

\subsection*{Acknowledgements}
 	The research of the first author was supported by the Austrian Science Fund (FWF) under grant P32832-N.


\appendix

\section{Proof of \cref{The : NCgen}}\label{sec:appendix}

\begin{proof}
 Let $x^* \in \widetilde{D}^\ast_{2} \Phi((\xb,0);(u,v))(y^*)$ and consider  
 $\{t_k\}_{k\in\N}\subset\R_+$, $\{u_k\}_{k\in\N},\{x_k^*\}_{k\in\N}\subset\mathbb X$, and $\{v_k\}_{k\in\N},\{y_k^*\}_{k\in\N}\subset\mathbb Y$
 with $t_k\searrow 0$, $u_k\to u$, $v_k\to v$, $x_k^*\to x^*$, $y_k^*\to y^*$, as well as
 \[(x_k^*,-y_k^*/\tau_k) \in \widehat{\mathcal N}_{\gph \Phi}(\xb + t_ku_k,t_k v_k)\]
 for all $k\in\N$ where we used $\tau_k := t_k \norm{u_k}$ for brevity of notation.
 \cref{lem:coderivatives_constraint_maps} yields 
 $x_k^*=\nabla g(\bar x+t_ku_k)^* y_k^*/\tau_k$ and $y_k^*\in\tau_k\widehat{\mathcal N}_D(g(\bar x+t_ku_k)-t_k v_k)$
 for each $k\in\N$. Taking the limit in $\tau_kx_k^*=\nabla g(\bar x+t_ku_k)^* y_k^*$, we find $y^*\in\ker\nabla g(\bar x)^*$.
 Combining this with a Taylor expansion and 
 denoting $\tilde w_k:=g(\bar x+t_ku_k)-t_k v_k$ gives us
\begin{subequations}
	\begin{align}
	\label{eqn2 : Domain}
		x_k^*  
		& 
		= 
		\nabla g(\xb)^* \frac{y_k^* - y^*}{\tau_k} + \nabla^2\langle y_k^*,g\rangle(\bar x)(u) + \oo(1),\\
	\label{eqn1 : Image}
		y_k^* 
		& 
		\in
		 \widehat{\mathcal N}_D(\tilde w_k)
		=
		 \widehat{\mathcal N}_D\left(
		 	g(\xb) + t_k \left(\nabla g(\xb) u - v + \oo(1)\right) 
		 	\right)
	\end{align}
\end{subequations}
for each $k\in\N$.

In the general case~\ref{item:general_estimate}, we readily obtain $y^* \in \mathcal N_{D}(g(\xb);\nabla g(\xb) u - v)$, i.e.,
\eqref{eq:multiplier_from_dir_lim_normal_cone_and_kernel}, as well as
\[
	x^* - \nabla^2\langle y^*,g\rangle(\bar x)(u) \in \Im \nabla g(\xb)^*
\]
due to the closedness of 
$\Im \nabla g(\xb)^*$.
Thus, \eqref{eq:2ordEstimNC} now holds for some $z^*\in[y^*]^\perp$ since $\mathbb Y = \spa y^* \oplus [y^*]^{\perp}$ 
is valid while $y^* \in \ker \nabla g(\xb)^*$ holds.

Let us now prove~\ref{item:general_estimate_+CQ}.
Using the notation from above, let us first assume that $\{z_k^*\}_{k\in\N}$, given by
$z_k^*:=(y_k^* - y^*)/\tau_k$ for each $k\in\N$, remains bounded.
Then we may pass to a subsequence (without relabeling)
so that it converges to some $z^*\in\mathbb Y$. We get
\[
	y^* + \tau_k z_k^* = y_k^* 
	\in  
	\widehat{\mathcal N}_D\left(g(\xb) + \tau_k (\nabla g(\xb)u - v + \oo(1))\right)
\]
and $z^* \in D\mathcal N_{D}(g(\xb),y^*)(\nabla g(\xb)u - v)$ follows.
Clearly, taking the limit in \eqref{eqn2 : Domain} yields \eqref{eq:2ordEstimNC} as well.

On the other hand, if $\{z_k^*\}_{k\in\N}$ does not remain bounded,
we pass to a subsequence (without relabeling) such that
$\tau_k/\norm{y_k^* - y^*} \to 0$ and $\hat{z}_k^* \to \hat{z}^*$
for some $\hat z^*\in\mathbb S_{\mathbb Y}$
where we used $\hat z_k^*:= (y_k^* - y^*)/\norm{y_k^* - y^*}$ for each $k\in\N$. 
Multiplying \eqref{eqn2 : Domain} by $\tau_k/\norm{y_k^* - y^*}$
and taking the limit yields $\nabla g(\xb)^* \hat{z}^* = 0$. 
Taking into account $(\tilde w_k - g(\xb))/\tau_k \to \nabla g(\xb)u - v$,
we get
\begin{equation}\label{eq:some_convergence_of_surrogate_sequences}
	\frac{\norm{\tilde w_k - g(\xb)}}{\norm{y_k^* - y^*}}
	=
	\frac{\norm{\tilde w_k - g(\xb)}}{\tau_k}\frac{\tau_k}{\norm{y_k^* - y^*}}
	\to
	0.
\end{equation}
Let us assume that $\nabla g(\bar x)u\neq v$.
Then, for sufficiently large $k\in\N$, we have $\tilde w_k\neq g(\bar x)$, so we can 
set $\hat q_k:=(\tilde w_k-g(\bar x))/\norm{\tilde w_k-g(\bar x)}$ for any
such $k\in\N$ and find $\hat q\in\mathbb S_{\mathbb Y}$ such that $\hat q_k\to\hat q$.
Moreover, we have
\[
	y^* + \norm{y_k^* - y^*} \hat{z}_k^* 
	= 
	y_k^* \in  \widehat{\mathcal N}_D\left(g(\xb) + \norm{\tilde w_k - g(\xb)} \hat{q}_k\right)
\]
from \eqref{eqn1 : Image}, so that \eqref{eq:some_convergence_of_surrogate_sequences}
yields $\hat{z}^* \in D_{\textup{sub}}\mathcal N_{D}(g(\xb),y^*)(\hat{q}) \subset D\mathcal N_{D}(g(\xb),y^*)(0)$.
This contradicts \eqref{eq:some:CQ}.
In case where $\nabla g(\bar x)u=v$ hold, \eqref{eq:some:CQ_2} is not applicable.
However, we still have
\[
	y^*+\norm{y_k^*-y^*}\hat z_k^*
	=
	y_k^* 
	\in 
	\widehat{\mathcal N}_D\left(g(\bar x)+\norm{y_k^*-y^*}\,\frac{\tilde w_k-g(\bar x)}{\norm{y_k^*-y^*}}\right),
\]
so that taking the limit $k\to\infty$ 
while respecting \eqref{eq:some_convergence_of_surrogate_sequences} 
yields $\hat z^*\in D\mathcal N_D(g(\bar x),y^*)(0)$
which contradicts \eqref{eq:some:CQ_1}.

Consider now the polyhedral case~\ref{item:polyhedral_estimate} and let us show that we are in a similar situation as in case~\ref{item:general_estimate_+CQ}
in a sense that we can always replace $(y_k^* - y^*)/\tau_k$ by a bounded sequence.
\Cref{lem:some_properties_of_polyhedral_sets}\,\ref{item:exactness_tangential_approximation} yields
the existence of a neighborhood $V\subset\R^m$ of $0$ such that
\[
	\mathcal T_D(g(\xb)) \cap V = \big(D - g(\xb) \big) \cap V,
\]
as well as the fact that $\mathcal T_D(g(\xb))$ is polyhedral.
Thus, from \eqref{eqn1 : Image} we conclude
\begin{align*}
	y_k^*
	& \in 
	\widehat{\mathcal N}_D\left(g(\xb) + t_k \left(\nabla g(\xb) u - v + \oo(1)\right)\right)\\
	& = 
	\widehat{\mathcal N}_{g(\xb) + \mathcal T_D(g(\xb))}\left(g(\xb) + t_k \left(\nabla g(\xb) u - v + \oo(1)\right)\right)\\
	& = 
	\widehat{\mathcal N}_{\mathcal T_D(g(\xb))}\left(\nabla g(\xb) u - v + \oo(1)\right)
	\ \subset \
	\mathcal N_{\mathcal T_D(g(\xb))}\left(\nabla g(\xb) u - v\right)
\end{align*}
for all sufficiently large $k\in\N$ and, particularly, $y^* \in \mathcal N_{\mathcal T_D(g(\xb))}\left(\nabla g(\xb) u - v\right)$ follows.
Moreover, $\mathcal N_{\mathcal T_D(g(\xb))}\left(\nabla g(\xb) u - v\right)$ is again a polyhedral set
and so it also enjoys the exactness of tangential approximation.
Thus, for sufficiently large $k\in\N$, we obtain $y_k^* -  y^* \in \mathcal T_{\mathcal N_{\mathcal T_D(g(\xb))}\left(\nabla g(\xb) u - v\right)}(y^*)$
which also gives us $(y_k^* -  y^*)/\tau_k \in \mathcal T_{\mathcal N_{\mathcal T_D(g(\xb))}\left(\nabla g(\xb) u - v\right)}(y^*)$.
Referring to \eqref{eqn2 : Domain} and taking into account that $\mathcal T_{\mathcal N_{\mathcal T_D(g(\xb))}\left(\nabla g(\xb) u - v\right)}(y^*)$ is also polyhedral,
we may now invoke Hoffman's lemma \cite[Lemma~3C.4]{DontchevRockafellar2014} to find some bounded sequence 
$\{z_k^*\}_{k\in\N} \subset \mathcal T_{\mathcal N_{\mathcal T_D(g(\xb))}\left(\nabla g(\xb) u - v\right)}(y^*)$
satisfying
 \[
 	\nabla g(\xb)^* z_k^* = x_k^* - \nabla^2\langle y^*_k,g\rangle(\bar x)(u) + \oo(1).
 \]
Thus, an accumulation point $z^*\in\R^m$ of $\{z_k^*\}_{k\in\N}$ satisfies  
 $z^* \in \mathcal T_{\mathcal N_{\mathcal T_D(g(\xb))}\left(\nabla g(\xb) u - v\right)}(y^*)$ and \eqref{eq:2ordEstimNC}.

Finally, consider the last case~\ref{item:polyhedral_estimate_order_2} and let $x^* \in D^\ast_{2} \Phi((\xb,0);(u,v))(y^*)$.
Similar arguments now yield \eqref{eqn2 : Domain} together with $y_k^* \in \widehat{\mathcal N}_D(w_k)$ for each $k\in\N$ where
\begin{align*}
	w_k
	& := 
	g(\bar x+t_ku_k)-t_k^2 v_k
	=
	g(\xb) + t_k \nabla g(\xb) u + t_k^2 z_k,\\
	z_k
	& := 
	\frac{\big( w_k - g(\xb) \big) / t_k - \nabla g(\xb) u}{t_k}
	=
	\nabla g(\xb) \frac{u_k - u}{t_k} + \frac12 \nabla^2 g(\xb)[u,u] -  v + \oo(1).
\end{align*}
Building upon the previous arguments from the polyhedral case and exploiting further the exactness of tangential approximation,
there also exists a neighborhood $W\subset\R^m$ of $0$ such that
\[
 	\mathbf T(u) \cap W 
 	= 
 	\mathcal T_{\mathcal T_D(g(\xb))}(\nabla g(\xb) u) \cap W 
 	=
 	\big(\mathcal T_D(g(\xb)) - \nabla g(\xb) u \big) \cap W
\]
 and $\mathbf T(u)$ is polyhedral.
 Consequently,  we have  $w_k - g(\xb) \in \mathcal T_D(g(\xb))$ and, hence, also
 $\big( w_k - g(\xb) \big) / t_k \in \mathcal T_D(g(\xb))$ for sufficiently large $k\in\N$.
 Similarly, we conclude that $z_k \in \mathbf T(u)$.
 Taking into account that for each cone $K$, $q \in K$, and $\alpha > 0$, one has $\mathcal T_K(q) = \mathcal T_K(\alpha q) $, 
 we find
 \begin{align*}
 	\mathcal T_D(w_k)
 	&=
 	\mathcal T_{g(\bar x)+\mathcal T_D(g(\bar x))}(w_k)
 	=
 	\mathcal T_{\mathcal T_D(g(\bar x))}((w_k-g(\bar x))/t_k)
 	\\
 	&=
 	\mathcal T_{\mathcal T_D(g(\bar x))}\left(\nabla g(\xb) u + t_k z_k \right)
 	=
 	\mathcal T_{\mathcal T_D(g(\bar x))-\nabla g(\bar x)u}\left(z_k\right)
 	=
 	\mathcal T_{\mathbf T(u)}\left(z_k\right)
 \end{align*}
 for all sufficiently large $k\in\N$, and we obtain $y_k^* \in \widehat{\mathcal N}_D(w_k) = \widehat{\mathcal N}_{\mathbf T(u)}(z_k)$.

 Since $\mathbf T(u)$ is polyhedral, so is $\gph \mathcal N_{\mathbf T(u)}$,
 see \cref{lem:normal_cone_map_of_polyhedral_set}, 
 and it can be written as a finite union of convex polyhedral sets, 
 say $C_1,\ldots,C_\ell\subset\R^m\times\R^m$.
 Thus, we have
 \[
 	(z_k,y_k^*) \in \gph \widehat{\mathcal N}_{\mathbf T(u)} \subset \gph \mathcal N_{\mathbf T(u)} = \bigcup_{j=1}^\ell C_j
 \]
 for sufficiently large $k\in\N$.
 We may pick an index $\bar j\in\{1,\ldots,\ell\}$ such that $(z_k,y_k^*) \in C_{\bar j}$ holds for infinitely many $k\in\N$
 and suppose that $C_{\bar j}$ can be represented as $C_{\bar j} = \{(z,y) \,|\, A z + B y \leq c\}$
 for some matrices $A$, $B$, as well as $c$ of appropriate dimensions.
 Hence, by passing to a subsequence (without relabeling), we get
 \[
 	A \nabla g(\xb)\frac{u_k - u}{t_k} \leq c - A \left(\frac12 \nabla^2 g(\xb)[u,u] - v + \oo(1)\right)- B y_k^*.
 \]
 For each $k\in\N$, a generalized version of Hoffman's lemma, 
 see \cite[Theorem~3]{Ioffe1979},
 now yields the existence of $s_k\in\mathbb X$ with
 \begin{align*}
 	A \nabla g(\xb) s_k &\leq c - A  ((1/2) \nabla^2 g(\xb)[u,u] - v+ \oo(1))- B y_k^*,\\
 	\norm{s_k} &\leq \beta\bigl\Vert c - A ((1/2) \nabla^2 g(\xb)[u,u] - v + \oo(1)) - B y_k^*\bigr\Vert
 \end{align*}
 for some constant $\beta > 0$ not depending on $k$.
 Thus $\{s_k\}_{k\in\N}$ is bounded and satisfies
 \[
 	\forall k\in\N\colon\quad
 	\bigl(\nabla g(\xb)s_k + 1/2 \nabla^2 g(\xb)[u,u] - v + \oo(1),y_k^*\bigr) \in C_{\bar j} \subset \gph \mathcal N_{\mathbf T(u)}.
 \]
 We may assume that $\{s_k\}_{k\in\N}$ converges to some $s\in\mathbb X$.
 Exploiting \eqref{eq:Tu_and_ws}, we infer
 \[
 	y^*_k \in \mathcal  N_{\mathbf T(u)}(w_s(u,v) + \oo(1)) \subset \mathcal N_{\mathbf T(u)}(w_s(u,v)),\qquad y^* \in \mathcal N_{\mathbf T(u)}(w_s(u,v))
 \]
 for all sufficiently large $k\in\N$ from polyhedrality of $\mathbf T(u)$ 
 and the definition of the limiting normal cone.
 
 The remainder of the proof now follows by invoking \cite[Theorem~3]{Ioffe1979} again 
 and using the same arguments like in the proof of case~\ref{item:polyhedral_estimate}.
%
%
%
\end{proof}

\end{document}